\documentclass[11pt]{article}

\usepackage[
    letterpaper,
    left=0.85in,
    right=0.85in,
    top=0.85in,
    bottom=1in
]{geometry}

\usepackage{graphicx}
\usepackage{tabularx}
\usepackage{threeparttable}
\usepackage{multirow}
\usepackage{url}
\usepackage{amsmath}
\usepackage{comment}
\usepackage{nameref}

\usepackage{placeins}
\usepackage{authblk}
\usepackage{hyperref}
\usepackage{float}
\usepackage{caption}
\usepackage{subcaption}
\usepackage{amssymb}
\usepackage{amsthm}
\usepackage[english]{babel}
\usepackage{tikz}
\usepackage{lipsum,lmodern}
\usepackage[most]{tcolorbox}
\usepackage{newtxtext,newtxmath}
\usepackage{color}
\usepackage{bm}
\usepackage{booktabs}
\usepackage{placeins}
\usepackage[nameinlink]{cleveref}

\title{\textbf{Learning Effective Models from Network Dynamics Data with Multiple Initial Conditions Using Weak Form SINDy}}

\author[1]{Moyi Tian$^*$}
\author[2]{Daniel A.~Messenger}
\author[1]{Vanja Dukic}
\author[1]{Nancy Rodr\'iguez}
\author[1]{David M.~Bortz}

\affil[1]{Department of Applied Mathematics, University of Colorado, Boulder, CO 80309 United States}
\affil[2]{Theoretical Division, Los Alamos National Laboratory, Los Alamos, NM 87545 United States}
\affil[*]{\normalfont Corresponding author: moyi.tian@colorado.edu}

\begin{document}

\maketitle

\begin{abstract}
Social systems consist of networks of individuals who influence one another through social interactions. Studying how processes evolve on these networks can help us better understand patterns of social behavior. We study a system that couples online and offline social activity and investigate how to learn effective models directly from data using Weak Form Sparse Identification of Nonlinear Dynamics (WSINDy), a method for discovering governing equations. We assess learning performance using data generated by a mean-field approximation model of a stochastic interaction process on networks and test how accurately the system can be recovered under different noise levels. Our results show that using more trajectories improves accuracy when noise is high, but only a small number of additional trajectories is needed to gain most of the benefit, with little improvement beyond that. We also learn effective ODE models from averaged stochastic data on networks. When traditional mean-field approximations fail, identifying continuum ODEs directly from stochastic processes yields efficient models that better match the data and provide deeper insight into the underlying dynamics.
\end{abstract}

\section{Introduction} \label{sec:intro}


The volume of online activities has boomed over the years, with about $73.2 \%$ of the global internet users, including $68.7 \%$ of the worldwide population, engaging in social networks in October 2025 \cite{Statista2025, Digital2026}. With the growing use of social media, studies increasingly show causal and reinforcing relationships between online and offline social activities. For example, studies have examined how images of the 2011 Egyptian revolution shared on social media contributed to the political movement \cite{Kharroub2016}. In the 2013 Kenyan and 2015 Nigerian presidential elections, grievances on Twitter were found to have ignited physical protests \cite{Mohler2020}. Moreover, the 2011 Occupy movement shows that online communication on Twitter and Facebook and in-person protests have bidirectional causality \cite{Bastos2015}. 

These real-world spillover effects highlight the need to study how social activity evolves in online engagement and offline physical actions, and how the two interact. Several studies have developed mathematical models to describe the spatiotemporal dynamics of large-scale social activities. Early work by Burbeck and collaborators introduced epidemiological modeling ideas to the study of riots \cite{Burbeck1978}. This line of research later gained renewed attention in the work of Bonnasse-Gahot and collaborators, who extended these ideas to a nonlocal model of the 2005 French riots to capture their spatiotemporal dynamics \cite{Bonnasse2018}. Several other models of social uprising, such as those that account for social tension \cite{Berestycki2010, Berestycki2015} or policing strategies \cite{Davies2013}, have also been developed and studied. While these earlier models offer various perspectives on social dynamics, they do not account for the increasing influence of social media. Motivated by recent work in epidemiology and opinion dynamics showing that network structure plays a crucial role in systems involving human interactions \cite{Kiss2017,Fibich2016,Noorazar2020,Porter2016,Peng2021}, \cite{Tian2025} proposed a two-layer online–offline compartmental framework. This model incorporates feedback between online information and offline physical interactions to capture coupled dynamics, with online interactions evolving on networks.

Despite advances in modeling the dynamics of social uprisings, a major challenge remains in linking mathematical models to data, which is necessary to uncover and validate the mechanisms governing real-world social behavior on complex networks. In \cite{Tian2026_LLM}, the authors took an initial step toward estimating model parameters for LLM engagement dynamics on networks using real-world data. In the related area of opinion dynamics, earlier work has developed probabilistic algorithms to infer agent-based models from opinion data \cite{Monti2020}, as well as nonparametric approaches to learn interaction kernels in mean-field models \cite{Lang2022}. There are also other recent works on predicting network dynamics from data, for example, \cite{Prasse2022} and \cite{Gao2022}. However, while these studies help connect data to models, they typically assume a fixed model structure or specific underlying mechanisms. 

More recently, data-driven methods such as dynamic mode decomposition (DMD) \cite{Proctor2018} and Sparse Identification of Nonlinear Dynamics (SINDy) \cite{Brunton2016_SINDy} have made it possible to learn effective models and governing equations directly from complex dynamical data. Related to networks, as an extension of SINDy, \cite{Basiri2025} proposed SINDy from Graph-structured data (SINDyG), which incorporates prior knowledge about the network, specifically the network's adjacency matrix, to guide model discovery. Another related work is \cite{Liu2025}, where the authors developed Group Similarity SINDy (GS-SINDy) and emphasized the need to incorporate multiple trajectory information for system learning.

However, many real-world datasets, particularly those arising from stochastic or network systems, are noisy, making derivative estimation unreliable and limiting the effectiveness of standard SINDy-type approaches. To address this challenge, methods based on weak formulations have been developed. Of particular interest is the Weak Form Sparse Identification of Nonlinear Dynamics (WSINDy), which was introduced in \cite{Messenger2021_WSINDy_ODE,Bortz2024}. This approach uses convolution with compactly supported test functions, allowing system identification without explicitly estimating derivatives. This approach has been remarkably successful and has led to many publications showing how the weak form can be used for both discovering equations and identifying systems. It has worked well in applications to ODEs \cite{Bortz2023_ODE_App,Chawla2025_ODE_App, Heitzmanbreen2025,Messenger2021_WSINDy_ODE, Messenger2024_ODE_App}, PDEs \cite{Lyons2025, Messenger2021_WSINDy_PDE, Messenger2022_OnlineWeakForm_PDE, Minor2025_PDE_App, Vasey2025}, SDEs \cite{Bortz2024, Minor2025_MeanFieldPDE}, reduced order modeling \cite{He2025, Tran2024, Wang2025_hydrodynamics}, and even closure discovery \cite{Messenger2025}. There is also recent work using WSINDy to learn mean-field equations for interacting particle systems \cite{Messenger2022_SDE}. A mathematical study on the asymptotic consistency of WSINDy was also recently published, providing a theoretical foundation for the method \cite{Messenger2024_WSINDy_ContinuumData}. For a good introduction to the methodology, we direct the interested reader to \cite{Messenger2024_SIAM_Intro} and the more technical introduction in \cite{Bortz2024}.

These developments motivate a closer examination of how data-driven system identification can be applied to models of coupled online–offline social dynamics. As an initial step, we examine a stochastic model and its mean-field approximation for a fully-mixed population, introduced in \cite{Tian2025}, and perform systematic computational studies to assess the ability of data-driven methods to recover the underlying governing dynamics from simulated data. We employ WSINDy for equation learning and use the fully-mixed population model as an entry point for assessing the feasibility of this task. This model with a ground-truth ODE system allows us to examine the impact of noise. We then extend the scope to networks, thereby increasing the structural complexity, which is closer to the real-world scenario where the underlying dynamics have no true differential equations, by experimenting with data from the stochastic model evolving on Erd\H{o}s-R\'{e}nyi random networks. Our results show that combining multiple trajectories with different initial conditions in model learning yields better learning performance with respect to error measurements, and surprisingly, adding more than two or three trajectories provides little additional benefit.

\textit{Outline:} In Section \ref{sec:model}, we introduce two models used for data-driven system identification of coupled online–offline dynamics: a stochastic model (Section \ref{subsec:sys_stoch}) and its fully-mixed continuum approximation (Section \ref{subsec:sys_fullymixed}). Section \ref{sec:wsindy} describes the system identification method, WSINDy. Results are presented in Section \ref{sec:results}, with Section \ref{subsec:results_ODE} focusing on learning noisy data generated from the continuum model and Section \ref{subsec:results_stoch} on that from the stochastic model. Finally, we conclude in Section \ref{sec:conclusion} with a discussion of implications and future directions.

\section{Models} \label{sec:model}
For completeness, we briefly discuss the model setup that is the focus of this work. Details of the motivation and analysis of this model can be found in \cite{Tian2025}. 

We study a coupled online-offline dynamics on population size $N$, in which each individual has a hybrid state that combines online and offline states. Online, individuals can be in one of three states: uninterested (\textbf{U}), engaged (\textbf{E}), or disengaged (\textbf{D}). Individuals in \textbf{U} can be influenced by their \textbf{E}-neighbors through online connections to become engaged, with a transmission rate $\tau$. When the population is fully mixed, meaning that everyone is connected to everyone else online, we consider a population-level transmission rate, $\beta = N \tau$. Each engaged individual can become disengaged at a rate $\gamma_i$, where the subscript $i$ stands for information layer. For the offline state, we consider the inactive state to be \textbf{UP}, the active state to be \textbf{P}, and the recovered absorbing state to be \textbf{R}. Similarly to the online setting, in the offline setting, individuals in state \textbf{P} can transition to state \textbf{R} at rate $\gamma_p$, where the subscript $p$ denotes the physical layer.

The two layers interact as follows: an individual with online state \textbf{E} can switch their offline state from $\bf UP$ to $\bf P$ at rate $\eta$. Conversely, higher offline activity ($\bf P$) increases the likelihood that people move from $\bf U$ to $\bf E$ online, with this feedback captured by the self-excitation rate $\theta$. Because this model focuses on online-to-offline spillover dynamics, it assumes that offline action is triggered solely by online engagement. Hence, the state in which an individual is uninterested online and participates offline is considered not to exist.

\begin{figure}[H]
    \centering

    \begin{subfigure}[t]{0.49\textwidth}
        \centering
        \includegraphics[width=\textwidth]{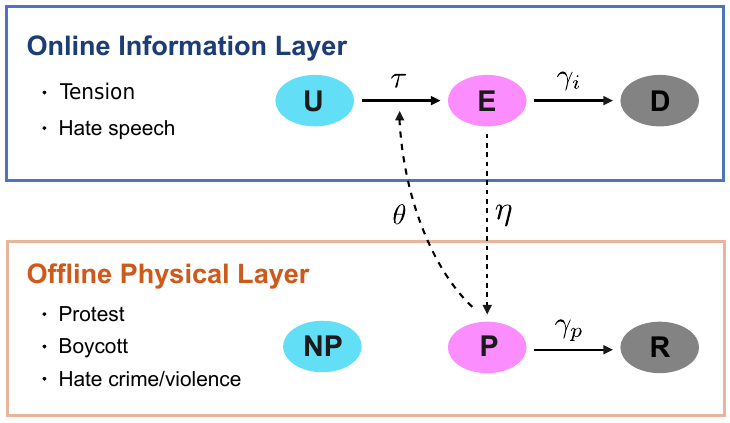}
        \caption{Model flow diagram}
        \label{subfig:flow_diagram}
    \end{subfigure}
    \hspace{20pt}
    \begin{subfigure}[t]{0.43\textwidth}
        \centering
        \includegraphics[width=\textwidth]{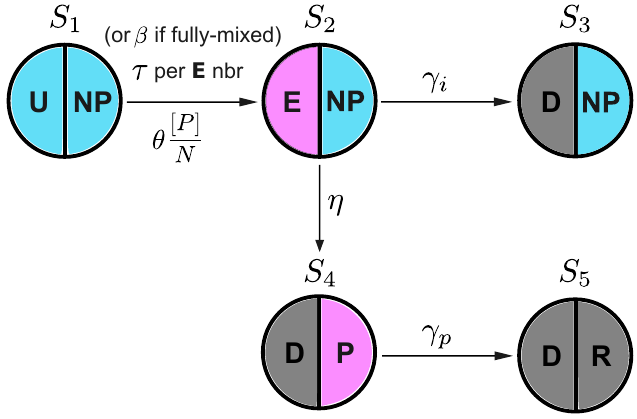}
        \caption{The Markov process}
        \label{subfig:markov_diagram}
    \end{subfigure}

    \caption{Illustrations of the model: (a) The flow within and between the online and offline spaces; (b) Diagram of the Markov stochastic process where $[P]$ denotes the number of people in offline $\mathbf{P}$ state and $N$ is the population size, so $[P]/N$ represents the offline engagement density. $\tau$ is the edge-level online transmission rate per contact, and $\beta = N \tau$, which is the population-level online transmission rate in the case of a fully-mixed population. $\eta$ is the online-to-offline transmission rate, $\theta$ is the self-excitement rate, and $\gamma_i, \gamma_p$ are the online and offline recovery rate, respectively.}
    \label{fig:model}
\end{figure}

We also assume the tension-inhibiting regime \cite{Bakhshi2021}, in which case once an individual's offline state becomes \textbf{P}, its online energy is automatically discharged, and the online state immediately switches from \textbf{E} to \textbf{D}. Finally, we assume permanent immunity in both the online and offline spaces, so that the \textbf{D} and \textbf{R} states are absorbing states. \Cref{table:1} summarizes the states and parameters in our system, and Figure \ref{subfig:flow_diagram} illustrates the flow diagram for this model.

\begin{table}[!ht]
\caption{Model states and parameters}
\label{table:1}
\centering
\small
\begin{tabular}{ll}
\toprule
\multicolumn{2}{c}{\textbf{Model states}} \\
\specialrule{1.2pt}{2pt}{2pt}

\multicolumn{2}{l}{\textbf{Online layer}} \\
\midrule
$\mathbf{U}$ & Uninterested individuals \\
$\mathbf{E}$ & Engaged individuals \\
$\mathbf{D}$ & Disengaged individuals (absorbing) \\[4pt]

\multicolumn{2}{l}{\textbf{Offline layer}} \\
\midrule
$\mathbf{UP}$ & Inactive individuals \\
$\mathbf{P}$ & Active individuals \\
$\mathbf{R}$ & Recovered individuals (absorbing) \\

\specialrule{1.2pt}{6pt}{2pt}

\multicolumn{2}{c}{\textbf{Model parameters}} \\
\specialrule{1.2pt}{2pt}{2pt}

\multicolumn{2}{l}{\textbf{Online layer}} \\
\midrule
Population-level transmission rate & $\beta$ \\
Edge-level transmission rate & $\tau$ \\
Recovery rate & $\gamma_i$ \\[4pt]

\multicolumn{2}{l}{\textbf{Offline layer}} \\
\midrule
Online-to-offline transmission rate & $\eta$ \\
Recovery rate & $\gamma_p$ \\[4pt]

\multicolumn{2}{l}{\textbf{Cross-layer interaction}} \\
\midrule
Self-excitation rate & $\theta$ \\
\bottomrule
\end{tabular}
\end{table}

\subsection{The stochastic model} \label{subsec:sys_stoch}

The process described above defines a continuous-time Markov chain with five hybrid states, denoted by $S_1, S_2, S_3, S_4, S_5$. State $S_1$ represents individuals who are uninterested in online activities and do not participate in offline activities. Individuals in $S_1$ may transition to $S_2$ by becoming engaged online while still not participating offline. From $S_2$, there are two possible transitions: to $S_3$, corresponding to individuals who disengage online and remain nonparticipants offline, or to $S_4$, corresponding to individuals who disengage online while beginning to participate offline. From $S_4$, individuals may transition to $S_5$, an absorbing state representing permanent online disengagement and done with participation offline. 

Throughout the paper, italic variables such as $U,E,D,P,R$ denote population fractions, while bold symbols such as $\mathbf{U},\mathbf{E},\mathbf{D},\mathbf{P},\mathbf{R}$ are used as labels for the corresponding individual-level states or compartments. Quantities in brackets follow standard count notation, where $[P]$ denotes the number of individuals in state $\mathbf{P}$, while $[UE]$ denotes the number of online edges between individuals in states $\mathbf{U}$ and $\mathbf{E}$.

We now formally define the stochastic process at the population level. Let $n_1(t),\ldots, n_5(t)$ be the counts of individuals in states $S_1,\ldots,S_5$ at time $t$, and $\mathbf{X}(t) = (n_1(t), \ldots, n_5(t))$ denote the vector of counts with $\sum_{i=1}^{5} n_i(t) = N$. According to the individual transition rules described above, there are
only four possible one-step transitions at the population level. Let $\mathbf{n}=(n_1,n_2,n_3,n_4,n_5)$ be the current population state, and let $q(\mathbf{n},\mathbf{n}')$ denote the transition rate from state $\mathbf{n}$ to $\mathbf{n}'$. Then we can explicitly present the transition rates, which are
\begin{align*}
    q\left((n_1,n_2,n_3,n_4,n_5),(n_1-1,n_2+1,n_3,n_4,n_5) \right) &= \tau [UE] + \theta \frac{n_1 n_4}{N},\\
    q\left((n_1,n_2,n_3,n_4,n_5),(n_1,n_2-1,n_3+1,n_4,n_5) \right) &= \gamma_i n_2,\\
    q\left((n_1,n_2,n_3,n_4,n_5),(n_1,n_2-1,n_3,n_4+1,n_5) \right) &= \eta n_2,\\
    q\left((n_1,n_2,n_3,n_4,n_5),(n_1,n_2,n_3,n_4-1,n_5+1) \right) &= \gamma_p n_4.
\end{align*}
All transition rates from $\mathbf{n}$ to states other than the four listed above are zero. The diagonal entry is then defined by
\[
q(\mathbf{n},\mathbf{n})
=
-
\sum_{\mathbf{n}'\neq \mathbf{n}}
q(\mathbf{n},\mathbf{n}').
\]
Since the rate of the first transition depends on $[UE]$, which is not determined by the count vector $\mathbf{X}(t)$ alone, the system is not closed at the level of population counts. Stochastic realizations of this process are generated using the Gillespie algorithm \cite{Kiss2017}, which is an exact simulation method for this continuous-time Markov chain. Figure~\ref{subfig:markov_diagram} illustrates the states and the transition probabilities between them.

\subsection{The fully-mixed population model} \label{subsec:sys_fullymixed}
Under the assumption that the population is fully mixed, which means ``random mixing" of the population, then over a small time interval, everyone has an equal chance of interacting with anyone else. In other words, interactions occur on a complete graph. For each population size $N$, denote 
\[\mathbf X_N(t)=\left(n_U(t),n_E(t),n_D(t),n_P(t),n_R(t)\right)\] 
as the population count process where $n_U = n_1, n_E = n_2, n_D = n_3+n_4+n_5, n_P = n_4$, and $n_R = n_5$ are the aggregated count variables, and denote
\[\bar{\mathbf X}_N(t)=\frac{1}{N}\mathbf X_N(t) = \left(U_N(t),E_N(t),D_N(t),P_N(t),R_N(t)\right)\]
as the corresponding density process. The fully-mixed stochastic model is considered under the standard density-dependent scaling $\tau_N=\frac{\beta}{N}$, where the edge-level transmission rate $\tau_N$ depends on $N$ and $\beta$ is fixed as $N\to\infty$. Equivalently, for each finite population size, $\beta=N\tau_N$.

In this special scenario, $[UE]$ is
$n_U n_E$. Hence the population-level rate for the count change from $\mathbf{U}$ to $\mathbf{E}$, for example, is
\[
    \tau_N n_U n_E + \theta \frac{n_U n_P}{N} 
    = \frac{\beta}{N} (N U_N)(N E_N) + \frac{\theta}{N} (N U_N) (N P_N)
    =
    N(\beta U_N E_N + \theta U_N P_N).
\]
The remaining transition rates can also be written as $N$ times functions of $\bar{\mathbf X}_N$. 

We denote $\mathbf x(t)=(U(t),E(t),D(t),P(t),R(t))$ as the solution of the following fully-mixed model by mean-field approximation, which is a deterministic system:
\begin{align}\label{model:fully_mixed}
\left\{\begin{array}{ll}
    \frac{dU(t)}{dt}&=-\beta UE - \theta UP,\\[6pt]
    \frac{dE(t)}{dt}& = \beta UE + \theta UP - (\eta+\gamma_i) E,\\[6pt]
    \frac{dD(t)}{dt}& = (\eta+\gamma_i) E,\\[6pt]
    \frac{dP(t)}{dt}& = \eta E-\gamma_p P,\\[6pt]
    \frac{dR(t)}{dt}& = \gamma_p P.\end{array}\right.
\end{align}
where $U, E, D, P, R \in [0,1]$ denote the aggregated population fractions associated with the corresponding compartments. The equations for $D$ and $R$ are included here to present the full population-level balance and the deterministic large-population limit. In the WSINDy learning experiments conducted in this work, however, we use only the active variables $(U,E,P)$ as training data, since $D$ and $R$ are absorbing variables determined by the transition flows and do not provide additional independent information for learning the dynamics.

According to \cite{Ethier_Kurtz1986} and \cite{Darling_Norris2008}, by the law of large numbers for density-dependent continuous-time Markov chains, assuming the initial density converges to $\mathbf x(0)$, the stochastic density process, $\bar{\mathbf X}_N(t)$, converges to $\mathbf x(t)$ uniformly on finite time intervals in probability as $N\to\infty$. Therefore, system~\eqref{model:fully_mixed} is the exact deterministic large-population limit of the fully-mixed stochastic process, although it is not an exact deterministic description for finite $N$. In the latter numerical experiments, we fix the population size $N$. For simplicity, when no confusion arises, we omit the subscript $N$ and write $U,E,D,P,R$ for the finite-population fractions.

This fully-mixed model, although largely simplified from the finer detailed stochastic process, is the simplest continuum ODE system that captures the key effects of the coupled dynamics and serves as an approximation that is accurate in the large-population limit, with deviations vanishing in probability on finite time intervals. Using analytical techniques of dynamical systems, solutions, stability, and the reproductive number can be studied, see \cite{Tian2025}. We use this model as a starting point to investigate how varying noise levels and the number of trajectories influence learning accuracy.

\section{Weak Sparse Identification of Nonlinear Dynamics (WSINDy)} \label{sec:wsindy}

Suppose that the true model to be learned is
\begin{equation*}
    \frac{d}{dt} \mathbf{x}(t) = \mathbf{F}(\mathbf{x}(t)), \quad \mathbf{x}(0) = x_0 \in \mathbb{R}^{1\times d}, \quad 0 \le t \le T,
\end{equation*}
where each element of $\mathbf{F}(\mathbf{x}(t))$ is a linear combination of (potentially nonlinear) candidate functions $\{ f_j \}_{j=1}^J$ with $J$ being the number of features, and the data to be used is $\mathbf{x} \in \mathbb{R}^{M \times d}$ with $M$ observed time points and $d$ states. To demonstrate the robustness of WSINDy, we learned the equations from artificial data generated by adding mean-zero Gaussian noise to the true solution. Specifically, for each observed state variable $X\in\{U,E,D,P,R\}$ and each sampled time point $t_m$, the noisy data are generated as
\[
    X_{\mathrm{noisy}}(t_m)=X_{\mathrm{true}}(t_m)+\varepsilon_m,
    \qquad
    \varepsilon_m \overset{\mathrm{iid}}{\sim}\mathcal{N}(0,\sigma_X^2),
\]
where the reported noise level $\alpha$ determines the noise standard deviation by
\[
    \sigma_X=\alpha\,\mathrm{RMS}(X_{\mathrm{true}})
    =
    \alpha\sqrt{\frac{1}{M}\sum_{m=1}^{M}X_{\mathrm{true}}(t_m)^2}.
\]
The noise variables are sampled independently across observed state variables and across time points. The noisy training data are not clipped to the interval $[0,1]$. The perturbed values are used directly in the WSINDy learning procedure, even if some observations fall slightly outside the physically meaningful range. This is to keep the additive noise model consistent with the weak-form integral formulation used in WSINDy, avoiding an additional nonlinear truncation in the training data. Thus a noise level of $0.1$ corresponds to additive Gaussian noise with standard deviation equal to $10\%$ of the RMS magnitude of the corresponding exact state solution.

Denote $\langle \cdot,\cdot \rangle$ as the $L^2$ inner product. For a family of test functions $\{ \phi_l \}_{l=1}^L$ where $L$ is the number of test functions, WSINDy constructs the weak-form linear system with $G \in \mathbb{R}^{L \times J}$, where $G_{l,j} = \langle \phi_l,f_j(\mathbf{X}) \rangle$, and $b \in \mathbb{R}^{L \times d}$, where $b_{l,i} = -\langle \phi_l',\mathbf{X}_i \rangle$. For each state equation $i$, WSINDy solves the generalized least squares 
\[
    G\mathbf{w}^{(i)} = \mathbf{b}^{(i)},
\]
where $\mathbf{b}^{(i)}$ is the $i$-th column of $b$ and $\mathbf{w}^{(i)}$ is the coefficient vector for the $i$-th state equation. Collecting these coefficient vectors gives the coefficient matrix
\[
    \mathbf{W}
    =
    \begin{bmatrix}
    \mathbf{w}^{(1)} & \mathbf{w}^{(2)} & \cdots & \mathbf{w}^{(d)}
    \end{bmatrix},
\]
so that the full system can be written as $G\mathbf{W}=b$. For regularization in least squares, we employ the modified sequential-thresholding least-squares algorithm (MSTLS) \cite{Messenger2021_WSINDy_PDE}, which iteratively adjusts the threshold for term selection based on data scales. 

For the test function, we use a piecewise polynomial of the form
\begin{equation} \label{eq:testfunc_poly}
    \phi(t,r) = 
    \begin{cases}
    C (r^2 - t^2)^p, &  t \in [-r, r] \\
    0, & \text{otherwise}.
    \end{cases}
\end{equation}
where $C=C(r,p)$ is chosen so that $\| \phi \|_2 = 1$. In the weak form, derivatives are transferred from the data to the smooth compactly supported test function, so the right-hand side uses $-\langle \phi_l', X_i\rangle$ rather than numerical derivatives of the observed data. For the polynomial test function in \eqref{eq:testfunc_poly}, the derivative $\phi'$ is computed analytically. For implementation, we set $p = 12$, which is chosen based on numerical experiments. Note that this choice aligns with the range of values implemented and suggested in recent work \cite{Lyons2025, Tran2025}. While optimal choice of the parameters used in test functions based on data remains an active frontier of research, \cite{Tran2025} provides an insight in hyperparameter selection for piecewise polynomial test functions, including the order. The idea behind is that increasing $p$ produces smoother test functions and consequently smoother Gram matrix $G_{l,j}$ and the right-hand side $b_{l,i}$, but overly large values amplify numerical errors due to the limitation of machine precision. This trade-off can also be interpreted in terms of the Maclaurin series, where a moderate choice of $p$ corresponds to a truncation that emphasizes dominant nonlinear effects while limiting sensitivity to higher-order terms.

To select the optimal support radius, we compute the Fourier transforms of the test function and the data, and pick the best match between the test function's spectrum and the data's spectrum. The centers of the test functions are placed at every point separated by $1/4$ of the radius to ensure overlap.

When using multiple trajectory data, we solve the least-squares problem:
\begin{equation*}
    \begin{bmatrix}
        G^1 & G^2 & \hdots & G^K
    \end{bmatrix}^T
    \mathbf{W} = 
    \begin{bmatrix}
        b^1 & b^2 & \hdots & b^K
    \end{bmatrix}^T
\end{equation*}
where $K$ is the number of trajectories included in the learning, and the same coefficient matrix $\mathbf{W}$ is used for all trajectories. This formulation gives equal weight to each trajectory, although trajectories generated from different initial conditions may differ in variance and in how strongly they inform the coefficient estimates. Explicitly modeling such differences is left for future work.

\section{Results} \label{sec:results}

In this section, we discuss the experiments and results of learning the system of ODEs from the data generated by the fully-mixed model \eqref{model:fully_mixed} and the stochastic model described in Section \ref{subsec:sys_stoch}. We study the learning performance of WSINDy across different numbers of trajectories for various initial conditions in the data. Motivated by practical considerations arising when online activities trigger offline events, we simulate data under the scenario of initial online engagement with no offline activity. We then apply WSINDy to learn the model from these data. The learning procedure uses a time series of the active variables $(U, E, P)$. In contrast, $D$ and $R$ are absorbing states, the inclusion of which does not provide additional information about the dynamics. Moreover, particularly in noisy data, we observe that incorporating these trajectories leads to less stable learned models, as amplified noise can increasingly contribute to an imbalance between inflow and outflow and violate mass conservation. Thus, the trajectory information from $\mathbf{D}$ and $\mathbf{R}$ is excluded from the learning. The regression does not directly enforce the full state constraints, such as $U+E+D=1$ or the consistency relations among $D$, $P$, and $R$.

We use the following error metrics to quantify the performance of the learning procedure.
\begin{enumerate}
    \item {\bf Equation Error.} 
    For each state equation $i$, we define the equation error as
    \[
        \frac{\|G\hat{\mathbf{w}}^{(i)}-\mathbf{b}^{(i)}\|_2}
        {\|\mathbf{b}^{(i)}\|_2},
    \]
    where $\hat{\mathbf{w}}^{(i)}$ is the estimated coefficient vector for the $i$-th state equation.
    
    \item {\bf Output Error.} 
    For each state $i$, we define the output error as
    \[
        \frac{\|\hat{\mathbf{X}}^{(i)}-\mathbf{X}^{(i)}\|_2}
        {\|\mathbf{X}^{(i)}\|_2},
    \]
    where $\hat{\mathbf{X}}^{(i)}$ is the inferred solution and $\mathbf{X}^{(i)}$ is the corresponding observed trajectory for the $i$-th state.
    
    Thus, these errors are reported state-wise rather than as a single norm over all state variables.
    
    \item If the ground-truth model and weights are available, then we also consider:
    \begin{enumerate}
        \item {\bf TPR (True Positive Ratio).} We use TPR introduced in \cite{Lagergren2020}: 
        \[ \text{TPR} = \frac{\text{TP}}{\text{TP}+\text{FN}+\text{FP}} \]
        where TP stands for the count of true positive terms learned, FN for the false negatives, and FP for the false positives. This score is equivalent to the Jaccard similarity between the learned support and the true support, and therefore measures the overlap between the nonzero terms selected by WSINDy and those appearing in the true model. True negatives are excluded, so the score is not dominated by correctly omitted zero terms when only a small subset of the candidate library appears in the true model.
        \item {\bf Parameter Error.} 
        We compute the parameter error as
        \[
            \frac{
            \left\|\operatorname{vec}\left(\hat{\mathbf{W}}\right)-\operatorname{vec}\left(\mathbf{W}\right)\right\|_2
            }{
            \left\|\operatorname{vec}\left(\mathbf{W}\right)\right\|_2
            },
        \]
        where $\hat{\mathbf{W}}$ and $\mathbf{W}$ denote the learned and true coefficient matrices across all state equations and candidate terms, respectively, and $\operatorname{vec}(\cdot)$ denotes vectorization using the same ordering of equations and terms.
    \end{enumerate}
\end{enumerate}

When generating data, either from the fully-mixed model or from the stochastic model, we sweep over initial conditions with small activation fractions in only one layer, in the set $I = \{0.01, 0.02, \ldots, 0.30\}$. The different numbers of trajectories refer to the number of distinct initial conditions used for training. For the fully-mixed model, each such trajectory is an ODE solution from one initial condition; for the stochastic model, each such trajectory is a sample-averaged trajectory obtained from stochastic Gillespie realizations. For example, when the initial activation is from online, the initial condition is $(U, E, P)(0) = (1-E(0), E(0), 0)$ where $E(0) \in I$. When the initial activation is from offline, the initial condition is $(U, E, P)(0) = (1-P(0), 0, P(0))$ where $P(0) \in I$ because our tension-inhibiting assumption of the model restricts that $D(0) = P(0)$ and thus by mass conservation, $U(0) = 1 - P(0)$. 

We select initial conditions using a progressive, nested maximin strategy. Starting from a single trajectory initialized at the center of the range, $0.15$, we iteratively add trajectories whose initial conditions maximize the minimum 
distance to those already selected. In particular, if $S_k$ denotes the current set of selected initial conditions, the next value is chosen as
\[
x_{k+1} = \text{argmax}_{x \in I}\min\nolimits_{y \in S_k} |x-y|.
\]
That is, at each step, we add the value whose distance to the closest previously selected point is maximal. This procedure produces a sequence of nested subsets that become increasingly uniform over the interval, while preserving all previously selected trajectories. We use this maximin rule as a space-filling heuristic, rather than as an information-optimal design for parameter estimation.

An additional implementation is to include constraints in the least squares to incorporate known mechanisms at the boundaries. For example, when $U=1$, we know $dU/dt = 0$, and when $E=0$, $dE/dt \ge 0$, etc. We test the inclusions of such constraints and found that although it helps enforce the learned system to remain within the physically meaningful bounds of $[0,1]$, the runtime becomes extremely long, making scaling up of the experiments unfeasible. Moreover, the output errors increase because of the additional restrictions applied, going against the goal of learning models that better match the data. Hence, we choose not to incorporate this implementation in our following report of the results. We leave further development along this line as a possible future direction.

\subsection{Fully-mixed model} \label{subsec:results_ODE}
We present results for varying levels of additive Gaussian noise. In this subsection, we report median values over repeated noise realizations for the fully-mixed noisy-data experiments. We use the median rather than the mean because the learned models can occasionally produce large errors under noisy data, leading to skewed performance distributions. Thus, the median provides a more robust summary of typical learning performance. When learning with WSINDy, we use a feature library that includes linear, bilinear, and quadratic terms. Including higher-order terms introduces rank deficiency in the regression, meaning some library terms are dependent on each other. We therefore include terms only up to quadratic to maintain identifiability among the terms with respect to the data. When generating inferred dynamics from the learned model, we restrict state values to be non-negative, which both helps avoid hanging in the forward solve and is a reasonable constraint given the physical interpretation of the solutions (ratios).

\begin{figure}[!htbp]
    \centering
    \begin{subfigure}[t]{0.3\textwidth}
        \centering
        \includegraphics[width=\textwidth]{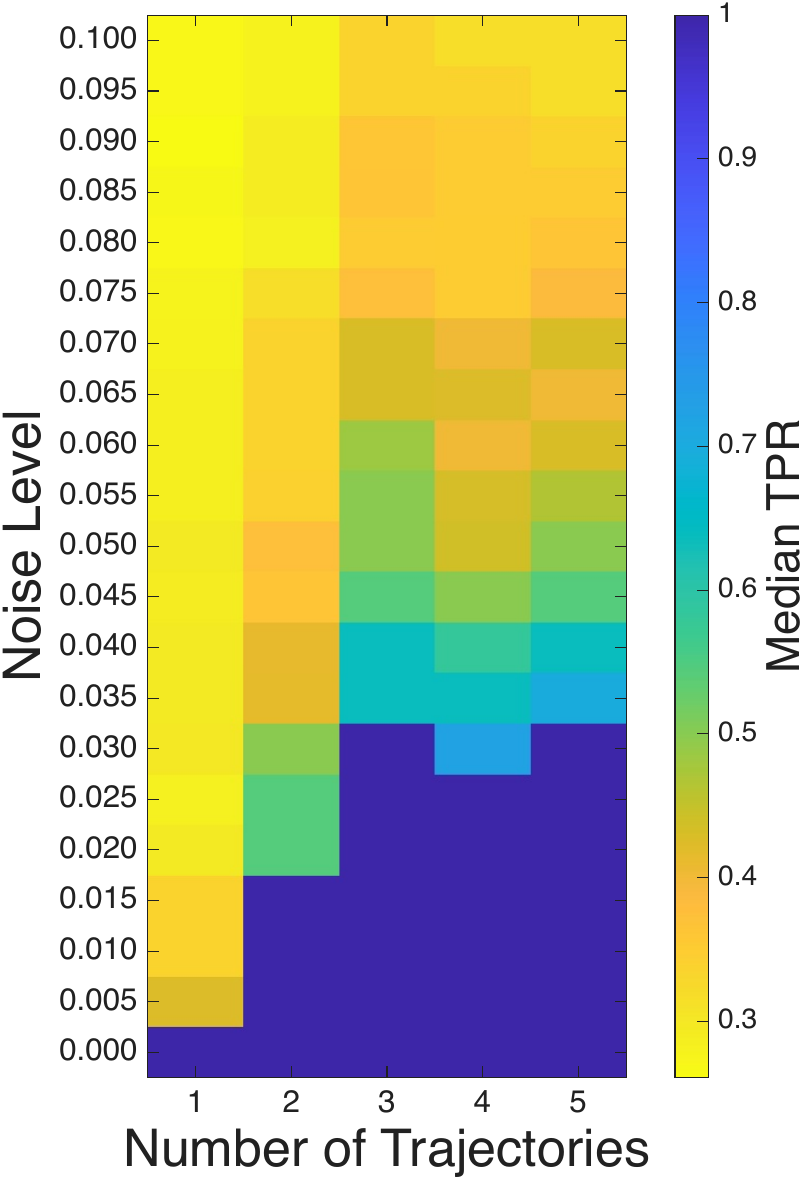}
        \caption{TPR}
        \label{subfig:E_TPS_Noise_NumTrajs}
    \end{subfigure}
    \hspace{60pt}
    \begin{subfigure}[t]{0.277\textwidth}
        \centering
        \includegraphics[width=\textwidth]{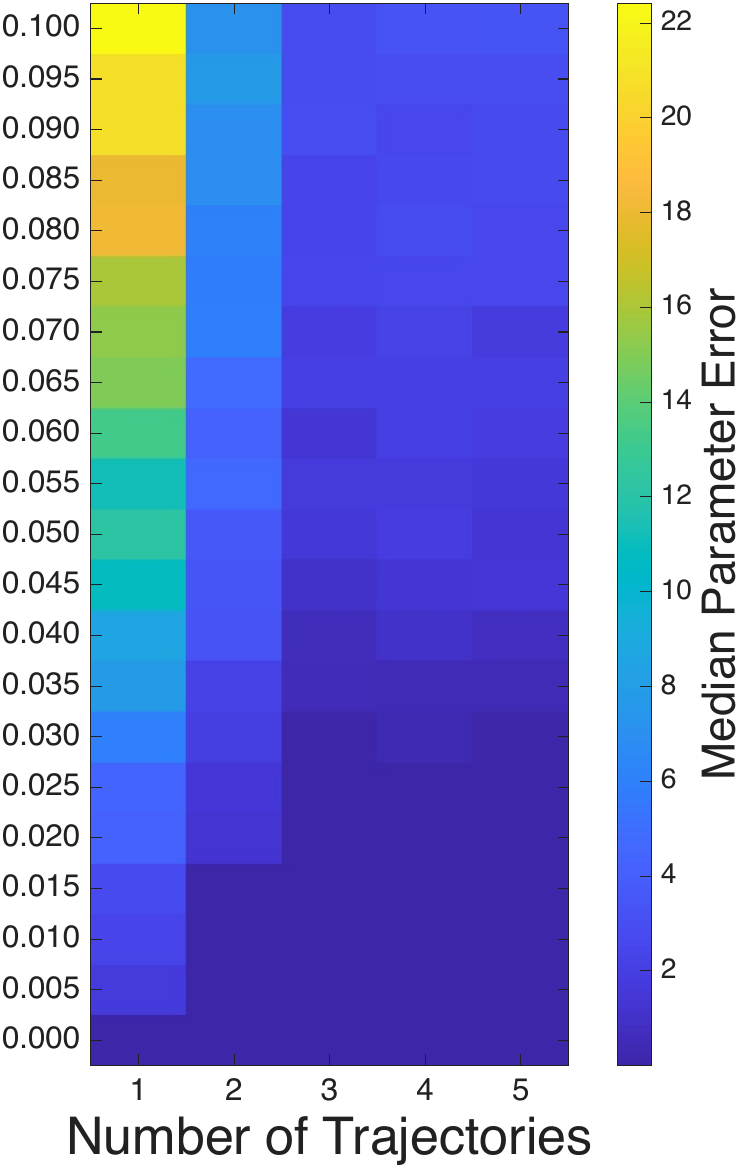}
        \caption{Parameter error}
        \label{subfig:E_ParamErr_Noise_NumTrajs}
    \end{subfigure}
    \caption{(a) Median TPR over $100$ realizations with varying initial conditions in $E(0)$. (b) Median parameter error over $100$ realizations with varying initial conditions in $E(0)$. The parameters used in the simulation are $\beta=0.5, \theta=0.4, \eta=0.2, \gamma_i=0.1, \gamma_p=0.3, t_{\max}=50, t_{\mathrm{num}}=1000$. When only $1$ trajectory is used for model learning, the accuracy of term selection and parameter inference drops significantly, even with minor noise. Providing $3$ trajectories substantially improves endurance under noise while maintaining high learning accuracy in the TPR and parameter error, but additional trajectories do not yield further improvement.}
    \label{fig:E_TPS_ParamErr_Noise_NumTrajs}
\end{figure}

Figure \ref{subfig:E_TPS_Noise_NumTrajs} presents the TPR across noise level and the number of trajectories. Consistent with the summary convention in this subsection, each grid is the median across results from $100$ noise realizations. The plot shows that model identification using a single trajectory performs well in the noiseless case but cannot withstand any noise. On the other hand, providing more trajectories improves noise handling, although after adding a small number of trajectories, no further improvements are observed with even more trajectories provided. Specifically, we observe that the accuracy of term recovery does not improve when the number of trajectories is increased to 3 in our computational experiments, which run up to 10 trajectories. To reduce redundancy and highlight the key takeaway, we only present results up to 5 trajectories in the figures. We observe the same takeaway when examining the parameter error, which is shown in Figure \ref{subfig:E_ParamErr_Noise_NumTrajs}. 

\begin{figure}[!htbp]
    \centering
    \begin{subfigure}[t]{0.33\textwidth}
        \centering
        \includegraphics[width=\textwidth]{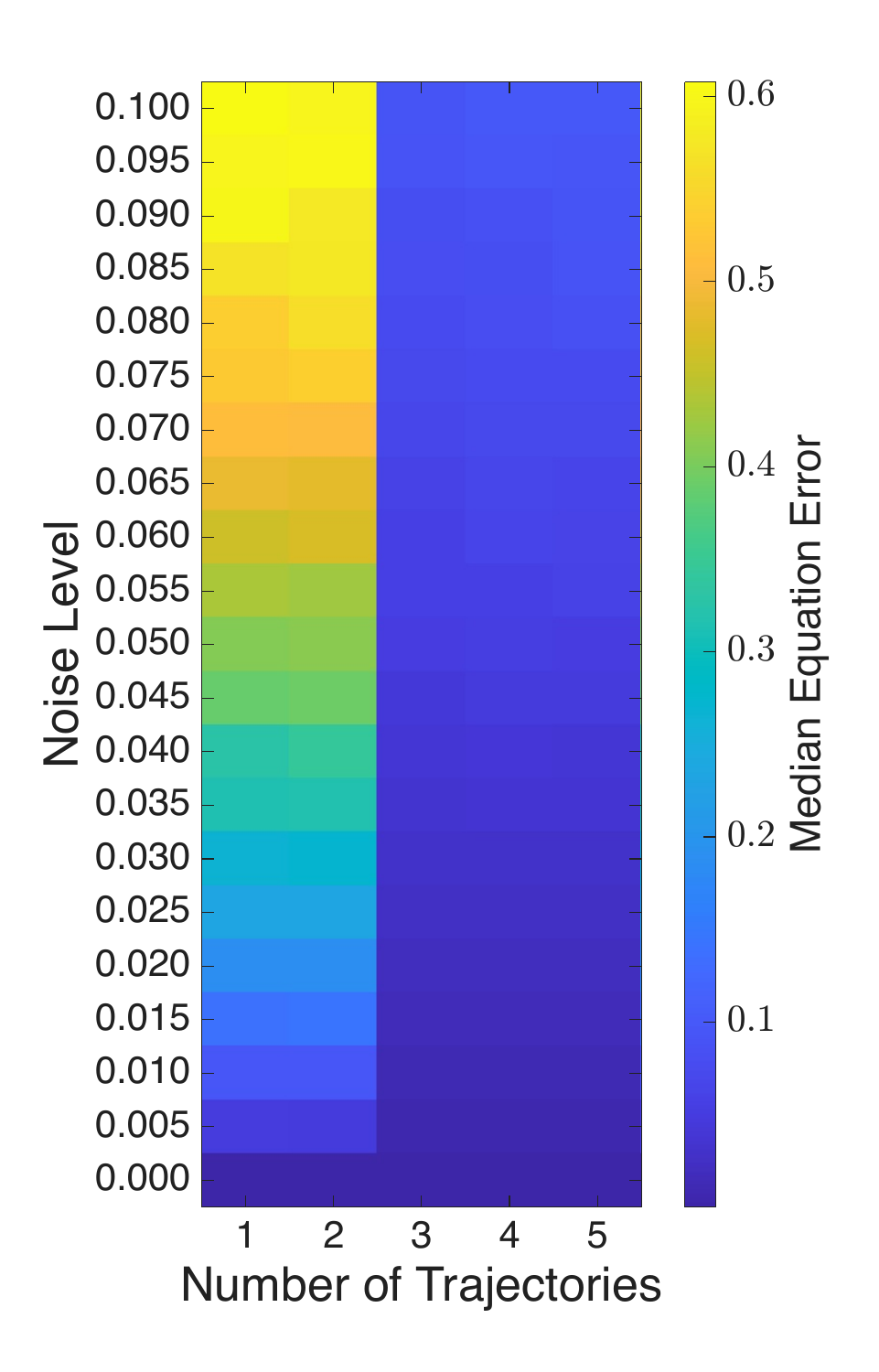}
        \caption{$E$ equation}
        \label{subfig:E_EquaErr_Noise_NumTrajs_E}
    \end{subfigure}
    \hspace{35pt}
    \begin{subfigure}[t]{0.33\textwidth}
        \centering
        \includegraphics[width=\textwidth]{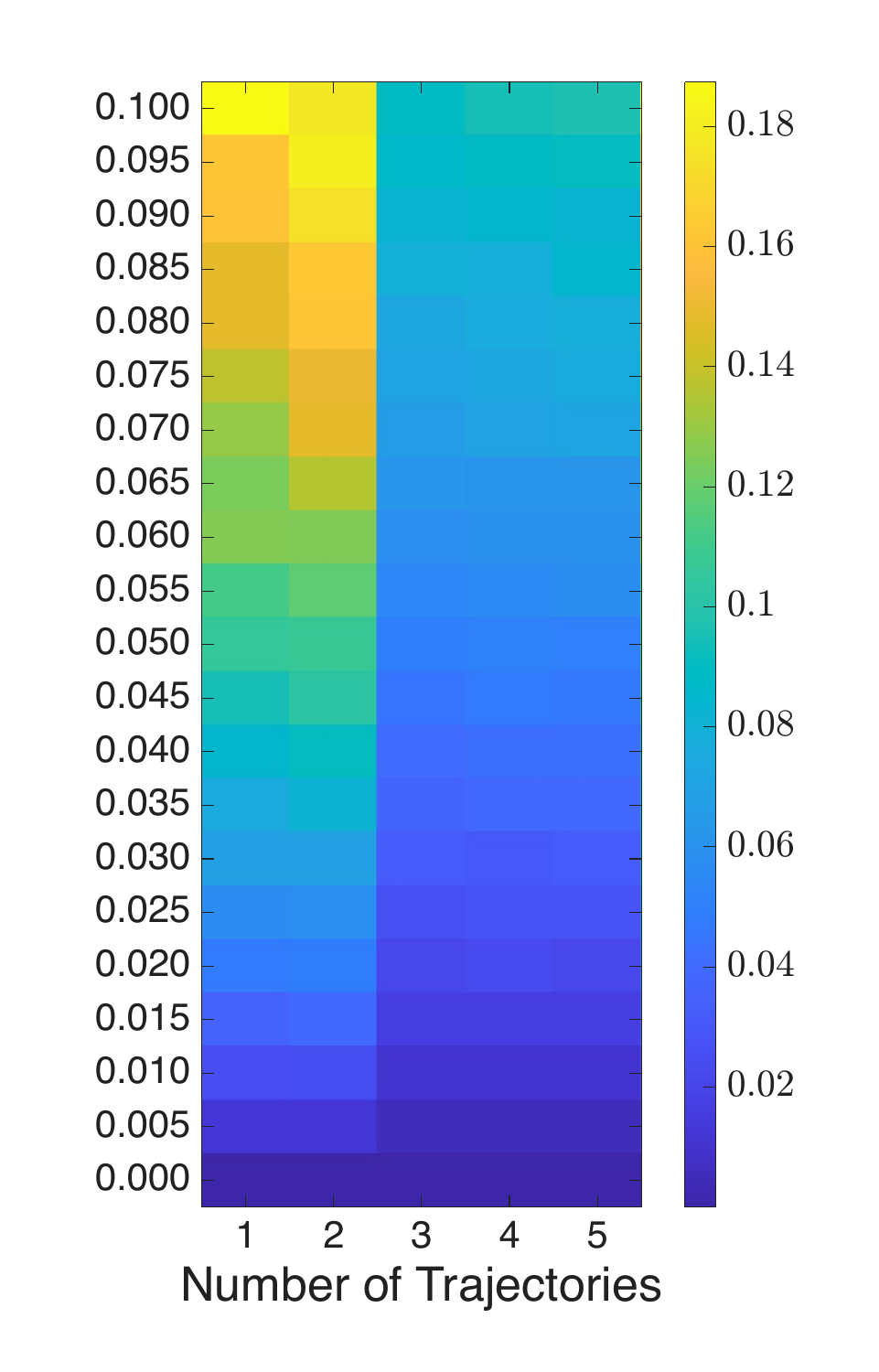}
        \caption{$P$ equation}
        \label{subfig:E_EquaErr_Noise_NumTrajs_P}
    \end{subfigure}
    \caption{Median equation error for (a) $E$ equation and (b) $P$ equation over $100$ realizations with varying initial conditions in $E(0)$. The same parameters are used for data generation. Significant improvement is observed when transitioning from $2$ to $3$ trajectories, with minor improvement beyond that.}
    \label{fig:E_EquaErr_Noise_NumTrajs}
\end{figure}

\begin{figure}[!htbp]
    \centering
    \begin{subfigure}[t]{0.33\textwidth}
        \centering
        \includegraphics[width=\textwidth]{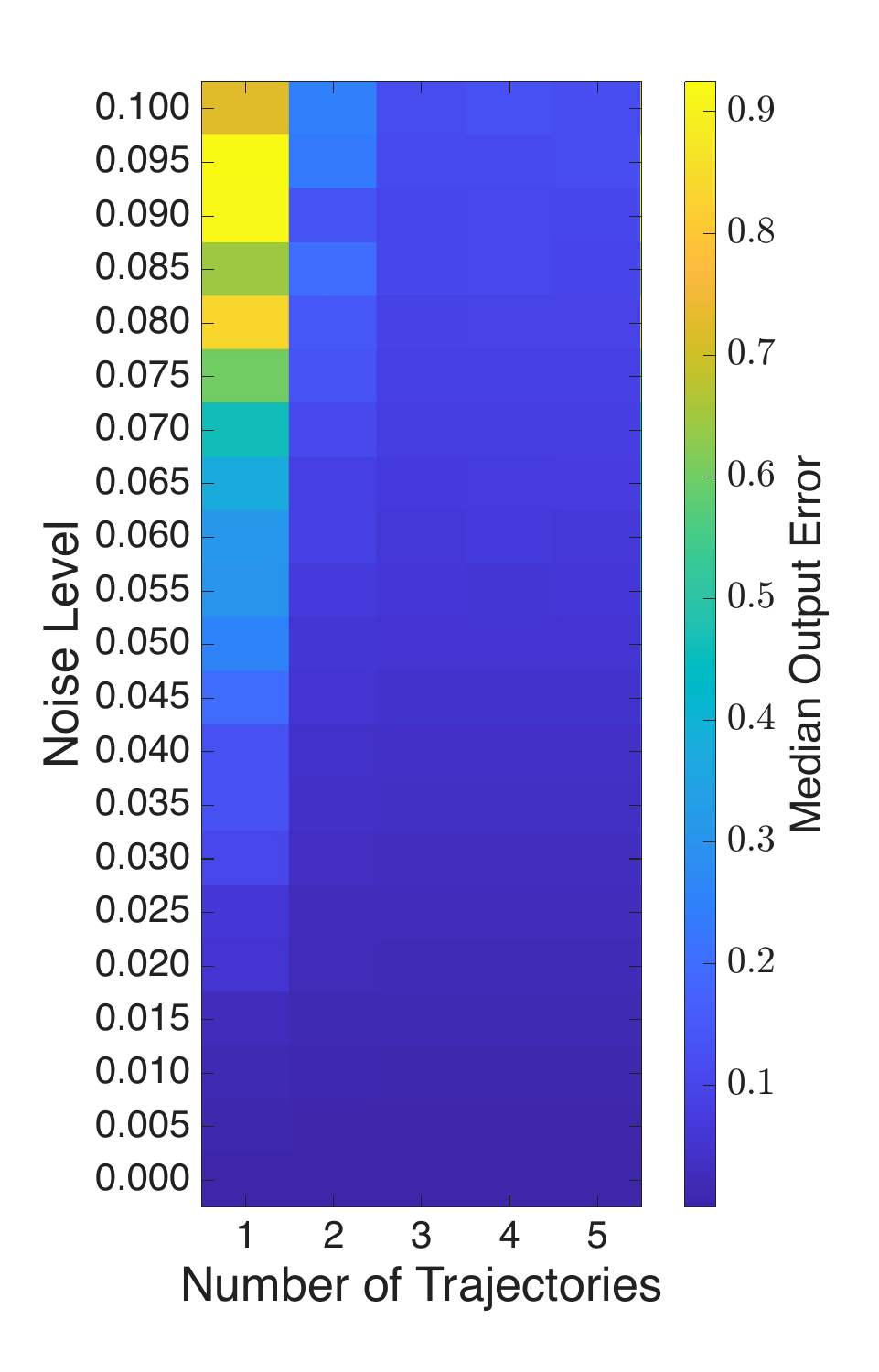}
        \caption{$E$ trajectory}
        \label{subfig:E_OutErr_Noise_NumTrajs_E}
    \end{subfigure}
    \hspace{35pt}
    \begin{subfigure}[t]{0.33\textwidth}
        \centering
        \includegraphics[width=\textwidth]{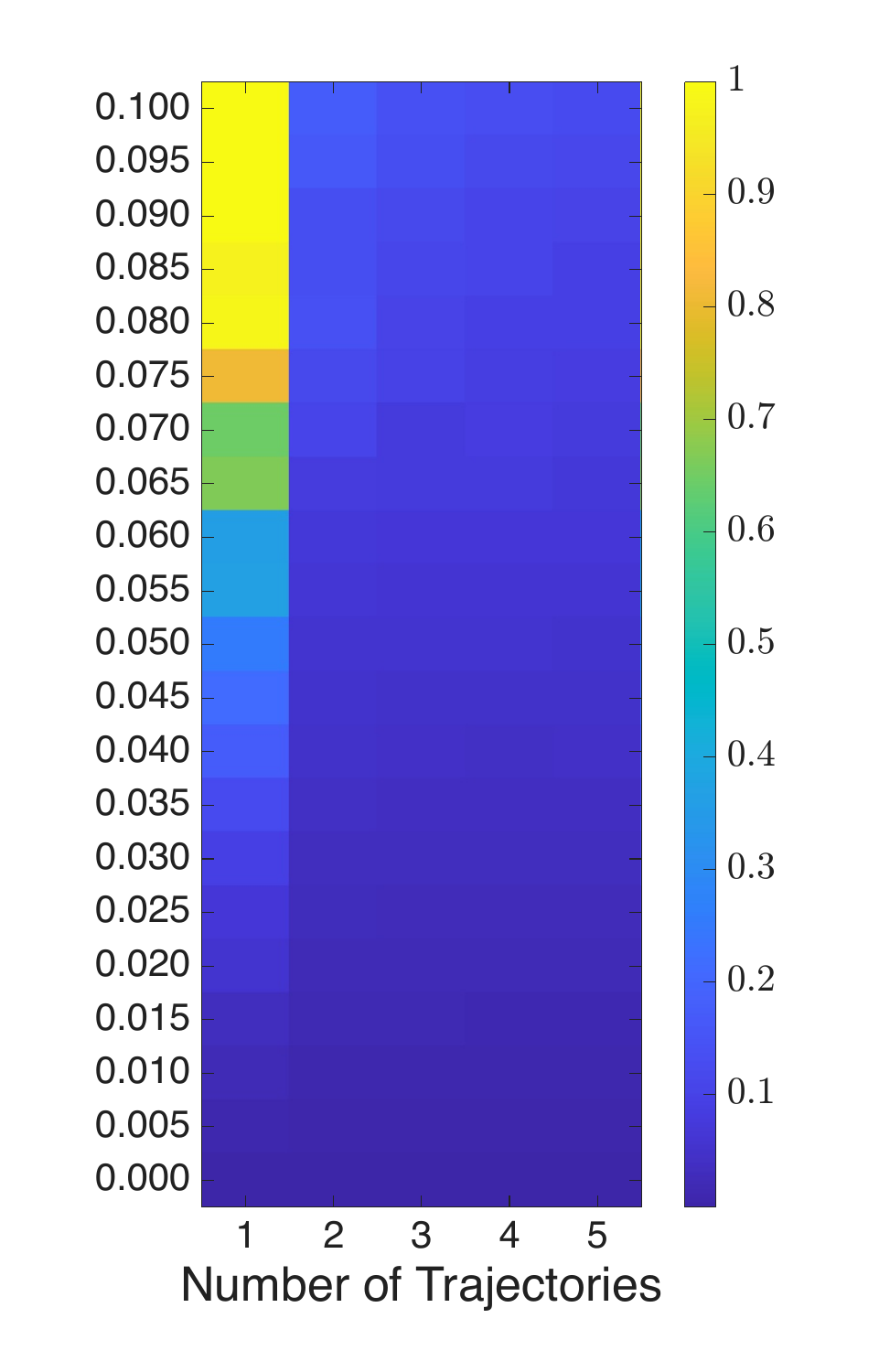}
        \caption{$P$ trajectory}
        \label{subfig:E_OutErr_Noise_NumTrajs_P}
    \end{subfigure}
    \caption{Median output error for (a) $E$ trajectory and (b) $P$ trajectory over $100$ realizations with varying initial conditions in $E(0)$. The same parameters are used for data generation. Improvements occur when the number of trajectories changes from $1$ to $2$ and from $2$ to $3$, but little change occurs beyond that.}
    \label{fig:E_OutErr_Noise_NumTrajs}
\end{figure}

We also examine and present the equation error and output error. These metrics are state-based, and we focus on the trajectories associated with states $\mathbf{E}$ and $\mathbf{P}$, which correspond to the quantities of interest for online and offline engagement, respectively. Figure \ref{fig:E_EquaErr_Noise_NumTrajs} shows the equation error for $E$ and $P$ equations. These plots reveal a clear improvement in performance as the number of trajectories increases from $2$ to $3$. Figure \ref{fig:E_OutErr_Noise_NumTrajs} presents the corresponding output error for $E$ and $P$ trajectories and shows two distinct transitions: first from $1$ to $2$ trajectories, and second from $2$ to $3$ trajectories. Together, these results indicate that using $3$ trajectories substantially improves learning performance, as measured by both equation and output errors. In contrast, adding more than $3$ trajectories yields little additional improvement in the median summaries. The corresponding variability at selected noise levels is shown in Appendix~\ref{app:variability_fullymixed}.

We conducted another set of experiments using generated data obtained by forward-solve with small fractions of offline activity and no initial online activity, namely, we set $E(0) = 0$ and vary $P(0)$. The results show a pattern similar to the above, in which $E(0)$ is small and $P(0)$ is zero. We thus omit the presentation of those results to reduce redundancy.

When the learned model is sufficiently accurate, it is natural to address identifiability. We discuss related considerations in Appendix \ref{app:eq_criterion}.

\subsection{Stochastic model} \label{subsec:results_stoch}

In this section, we present the learning results for dynamics generated by the stochastic model on random networks. Specifically, we investigate Erd\H{o}s-R\'{e}nyi random networks with $N=1000$. Each stochastic training trajectory is obtained by averaging $100$ Gillespie realizations generated from the same initial condition and network. We use the time interval $[0,t_{\mathrm{end}}]$ with $t_{\mathrm{end}}=40$ for learning and testing, and only Gillespie realizations that reach time $t_{\mathrm{end}}$ are retained so that every sample-averaged trajectory is defined over the full interval. The choice of $100$ realizations was made empirically, as preliminary tests indicated that this number was sufficient to stabilize the sample mean while keeping the simulation runtime reasonable. Therefore, the target of learning in this subsection is an effective deterministic model for the sample-averaged dynamics, rather than an individual sample path of the continuous-time Markov chain. 

Due to the stochastic nature of the model, we do not add additional noise to the data. For each WSINDy fit, the multiple training trajectories are generated from the same network realization but different initial conditions. Different network realizations are treated as separate replicates and lead to a distribution of learned models. Unlike the fully-mixed noisy-data experiments, the variability here comes from stochastic simulation and random network realizations rather than from additive observational noise scaled to the magnitude of each trajectory. We therefore report mean errors in this subsection to summarize average performance over the network and test-data replicates. To assess how well the learned model represents the data, we use output error as a metric, as a practical goal is to recover and predict the dynamics accurately. In WSINDy, we use a feature library that includes linear, quadratic, and cubic monomials to capture the additional complexity introduced by the network structure. As previously mentioned, libraries with higher-order terms exhibit rank deficiency issues and are therefore not suitable for our data.

We consider data generated on Erd\H{o}s-R\'{e}nyi random networks with varying average degree, $k$. Specifically, we use the $G(N,M)$ model, where the number of edges $M$ is chosen as the nearest integer to $Nk/2$. Thus, the realized average degree $2M/N$ is approximately $k$. We comment that this construction is closely related to the $G(N,p)$ model with $p = k/(N-1)$, for which the degree distribution is approximately Poisson with mean $k$ when $N$ is large and $k$ is small relative to $N$. As real networks are usually sparse, with a small average degree relative to their size \cite{Newman2003}, we focus on small average degree $k$ in the following results.

\begin{figure}[!htbp]
    \centering
    \includegraphics[width=0.53\linewidth]{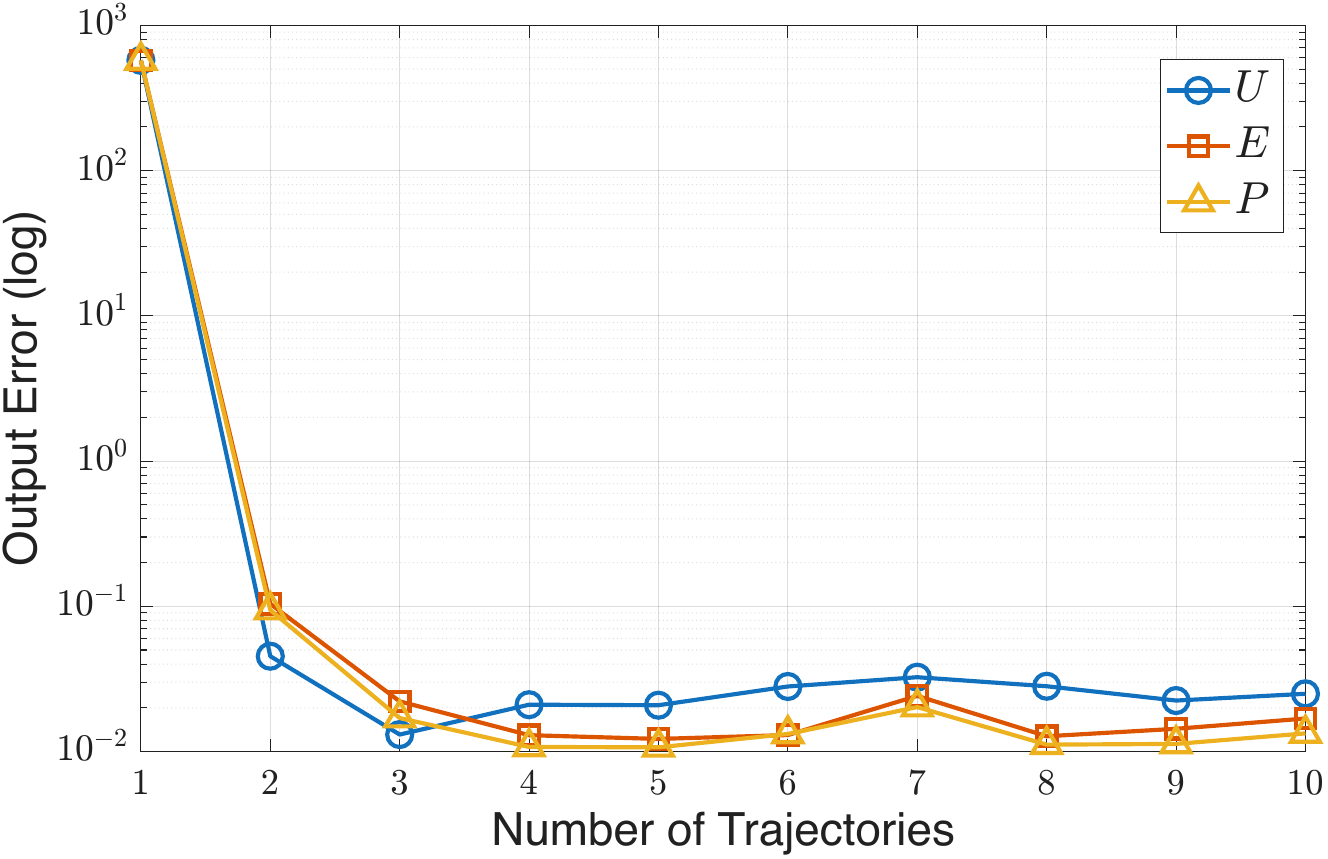}
    \caption{Mean output error relative to training data in base-$10$ logarithmic scale over number of trajectories for Gillespie simulations on a fixed Erd\H{o}s-R\'{e}nyi network with $N=1000, k=5, \tau=0.2, \theta=0.5, \eta=0.4, \gamma_i=0.1, \gamma_p=0.2, t_{\mathrm{end}}=40$ using library with polynomials from degree $1$ to $3$. Increasing the number of trajectories from 1 to 2 demonstrates significant improvement.}
    \label{fig:ER_TrajErr}
\end{figure}

Figure~\ref{fig:ER_TrajErr} shows one test example for data generated on one fixed Erd\H{o}s-R\'{e}nyi network with $N = 1000$ and $k = 5$. It illustrates the mean output error relative to the training data on a log scale for each state across different numbers of trajectories used in the learning. The result shows a substantial drop in error from $1$ to $2$ trajectories, but subsequent improvements are smaller, with fluctuations due to the stochastic nature.

For many models trained with a minimal number of trajectories, $1$ or $2$, the inferred solutions often deviate significantly from the real data. We observe that the learned model is sometimes too stiff for solvers to handle, or that inferred solutions, even when successfully solved, may go outside the $[0,1]$ bounds and thus lose physical meaning. Alternatively, inferred solutions may remain constant, far from expected behavior. We consider all of these scenarios to be failures of the learning process. For diagnostic purposes, we separate these failures into three disjoint categories: integration failures, including solver exceptions, nonfinite outputs, or incomplete time grids; out-of-range solutions, excluding cases already counted as integration failures; and flat solutions, excluding cases already counted in the first two categories. The failure rate reported below is the union of these three failure modes. To quantify how often such failures occur across different network densities and numbers of trajectories used for learning, we record the ratio of failed runs to total runs, which we call the failure rate. In these runs, for each average degree $k$, we randomly generate $100$ Erd\H{o}s-R\'{e}nyi networks where we sweep learning each across different numbers of trajectories. We then compare the inferred solutions against $10$ test sample-averaged trajectories, selected from trajectories with initial conditions not used in the training process, for each learned model and network. 

\begin{figure}[!htbp]
    \centering
    \begin{subfigure}[t]{0.34\textwidth}
        \centering
        \includegraphics[width=\textwidth]{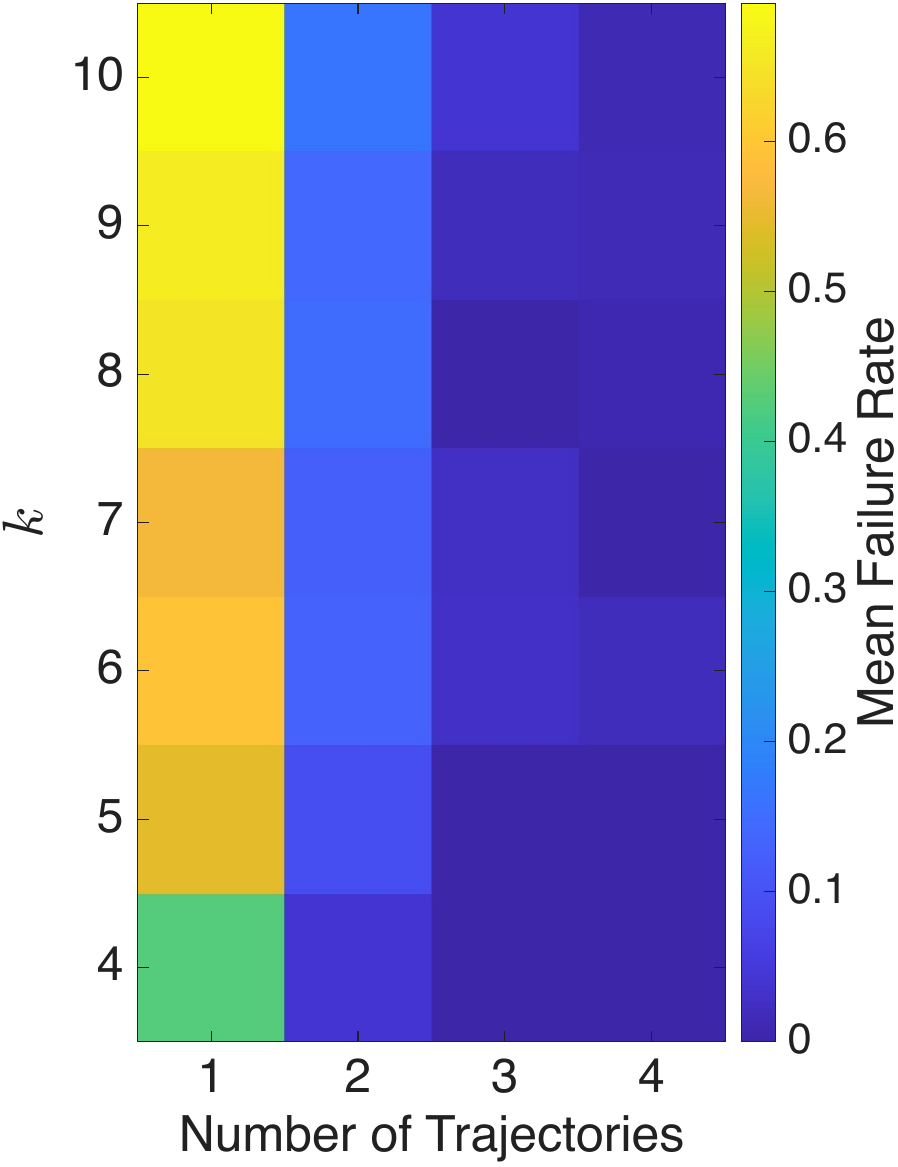}
        \caption{Failure rate heatmap}
        \label{subfig:ER_FailRate_left}
    \end{subfigure}
    \hspace{40pt}
    \begin{subfigure}[t]{0.38\textwidth}
        \centering
        \includegraphics[width=\textwidth]{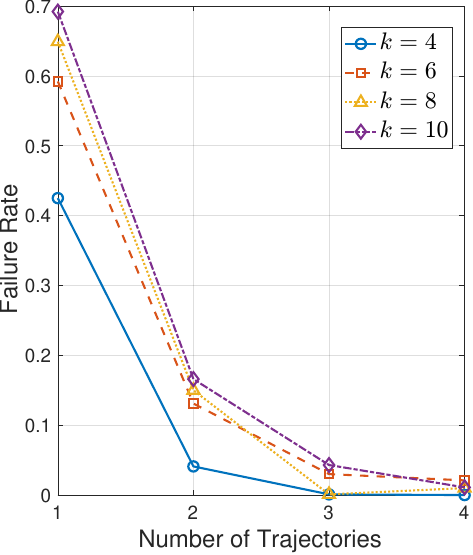}
        \caption{Failure rate for different $k$}
        \label{subfig:ER_FailRate_right}
    \end{subfigure}
    \caption{Failure rate heatmap (a) and example line plots for different $k$ (b). The same parameters are used for data generation as presented in Figure~\ref{fig:ER_TrajErr}. A significant failure rate is observed when only $1$ trajectory is used for learning, with the percentage of failures increasing with $k$. The failure rate drops significantly when the number of trajectories is increased to $2$, and again when it is further increased to $3$. After increasing the training set to $3$ trajectories, only small additional improvement is observed.}
    \label{fig:ER_FailRate}
\end{figure}

Figure~\ref{subfig:ER_FailRate_left} shows the heatmap of failure rate over $k$ and number of trajectories. Each grid shows an average over $1000$. Figure~\ref{subfig:ER_FailRate_right} presents several line plots of the failure rate over the number of trajectories corresponding to a few different $k$ values. The results show the first significant drop in the failure rate when increasing the number of trajectories from $1$ to $2$, and the second drop when increasing it from $2$ to $3$. Averaged over $k$, the rates of integration failures, bound violations, and flat solutions are respectively $34.4\%$, $16.3\%$, and $8.4\%$ of all test runs when $1$ trajectory is used; $6.2\%$, $1.9\%$, and $4.0\%$ when $2$ trajectories are used; and $0.4\%$, $0.01\%$, and $1.4\%$ when $3$ trajectories are used. For four trajectories, these rates further decrease to $0.03\%$, $0.1\%$, and $0.7\%$, respectively. Thus, the reduction from $1$ to $3$ trajectories is primarily driven by the disappearance of integration failures and bound violations, while the small residual failure rate for $3$ or more trajectories is mostly due to flat inferred solutions. Only small additional changes occur after increasing the training set to $3$ trajectories.

\begin{figure}[!htbp]
    \centering
    \begin{subfigure}[t]{0.333\textwidth}
        \centering
        \includegraphics[width=\textwidth]{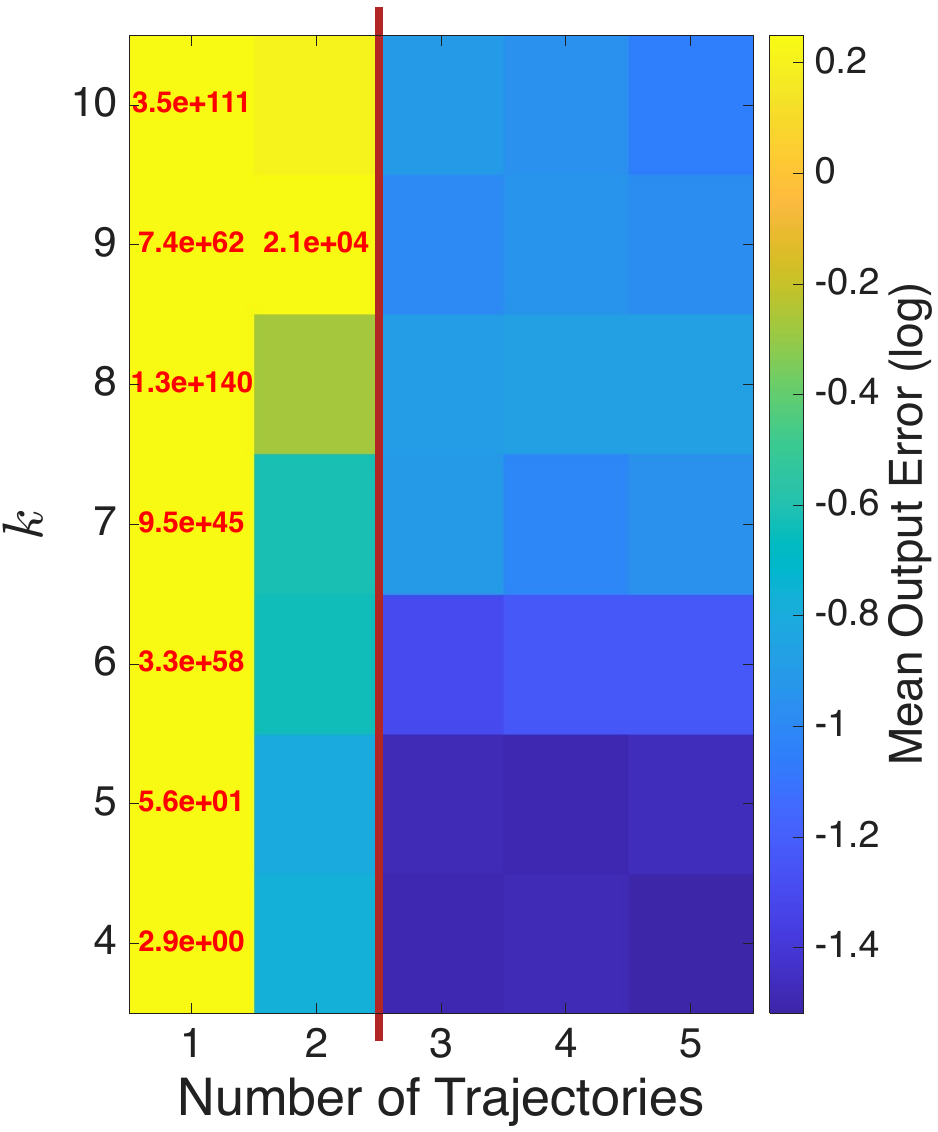}
        \caption{$U$ trajectory}
        \label{subfig:ER_Traj_Err_U}
    \end{subfigure}
    \hfill
    \begin{subfigure}[t]{0.295\textwidth}
        \centering
        \includegraphics[width=\textwidth]{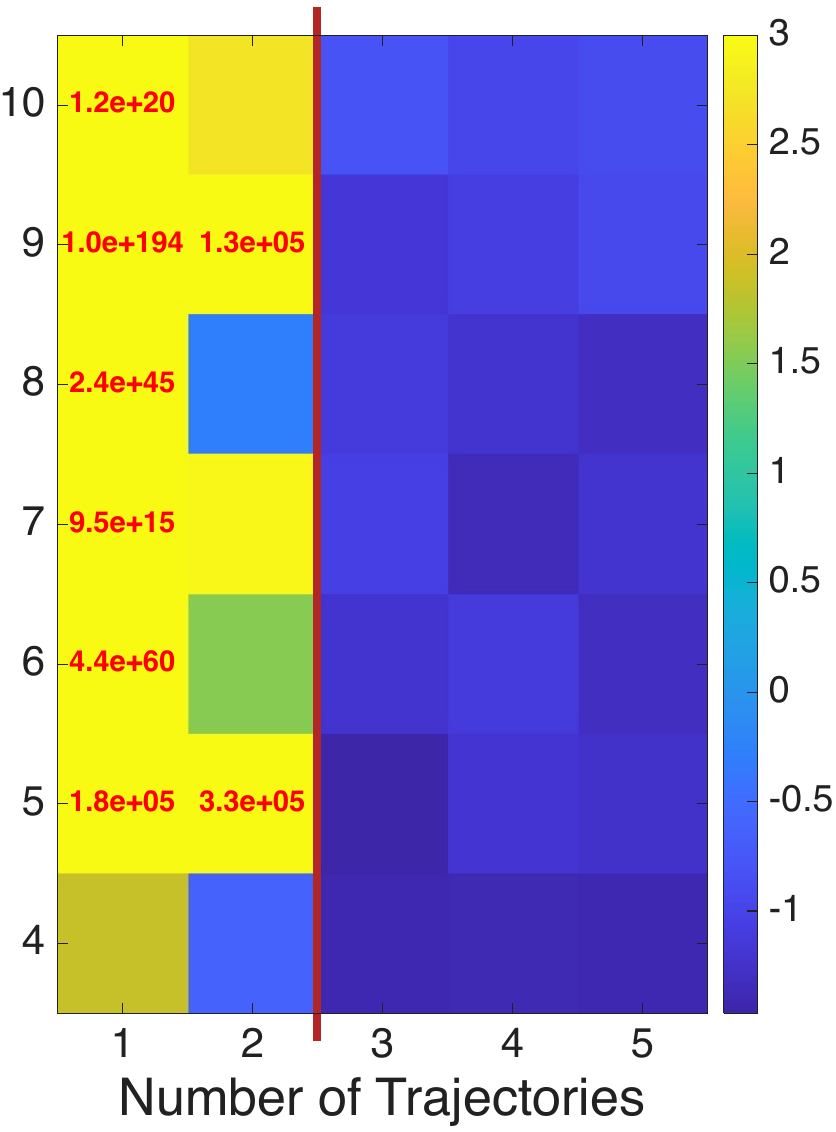}
        \caption{$E$ trajectory}
        \label{subfig:ER_Traj_Err_E}
    \end{subfigure}
    \hfill
    \begin{subfigure}[t]{0.295\textwidth}
        \centering
        \includegraphics[width=\textwidth]{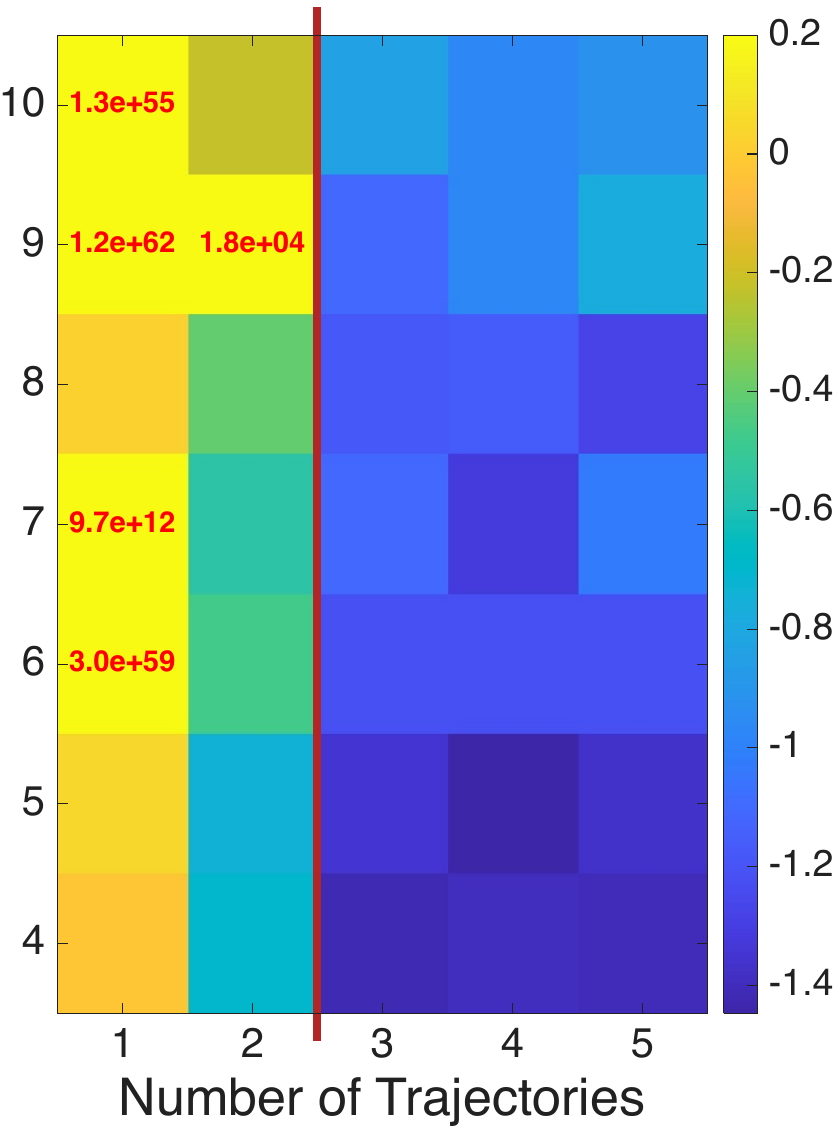}
        \caption{$P$ trajectory}
        \label{subfig:ER_Traj_Err_P}
    \end{subfigure}
    \caption{Mean output error for (a) $U$ trajectory, (b) $E$ trajectory, and (c) $P$ trajectory across average degree $k$ and number of trajectories on Erd\H{o}s-R\'{e}nyi random networks. Each grid shows an average over $1000$ values, with $100$ realizations including a random network generation each time, and each is then computed against $10$ test datasets. The same parameters are used for data generation as above. We employ colorbar cutoffs to avoid undue influence from extreme values. Across all three states, we observe an improvement when increasing from $2$ to $3$ trajectories, while changes beyond $3$ trajectories are small and inconsistent.}
    \label{fig:ER_Traj_Err}
\end{figure}

We present in Figure \ref{fig:ER_Traj_Err} the mean output error with respect to each state for $100$ independent network realizations, computed against $10$ test data sets, so each grid entry is an average over $1000$ values. The results provide additional evidence that improvement is observed when the number of training trajectories is increased to $3$, with only small and inconsistent changes beyond that. Representative cross-sectional boxplots at $k=5$ for these heatmaps are provided in Appendix~\ref{app:variability_stochastic}, showing the variability across network realizations and test sample-averaged trajectories behind the mean summaries.

Note that in Figures \ref{fig:ER_FailRate} and \ref{fig:ER_Traj_Err}, we show results on refined regions from $k=4$ to $10$ and number of trajectories from $1$ to $4$ or $5$ in order to highlight where the key transitions happen. In these experiments, we have run over the range $k = 4$ to $20$ and the number of trajectories from $1$ to $10$. The qualitative behaviors beyond the presented regimes exhibit the same trend and are therefore omitted from these figures for the sake of space.

There is also a reason related to the accuracy of mean-field approximations for focusing on sparse networks in the following comparison. For denser networks, mean-field approximations have been proposed and studied as a promising approach for deriving relatively accurate models to predict population-level dynamics of stochastic processes on heterogeneous networks \cite{Pastor-Satorras2001, Eames2002}. However, as networks become sparse, mean-field models may become less accurate because they neglect the differentiation of the individual nodes and their local connectivity, which become increasingly important in networks with low average degree \cite{Kiss2017}. We therefore focus our discussion on the sparse-network regime, where these local network effects are more pronounced, highlighting the data-driven modeling as a useful alternative. Related details on how network density affects the accuracy of mean-field approximations, specific to our system of studies, are discussed in \cite{Tian2025}.

\begin{figure}[H]
    \centering
    \begin{subfigure}[t]{0.315\textwidth}
        \centering
        \includegraphics[width=\textwidth]{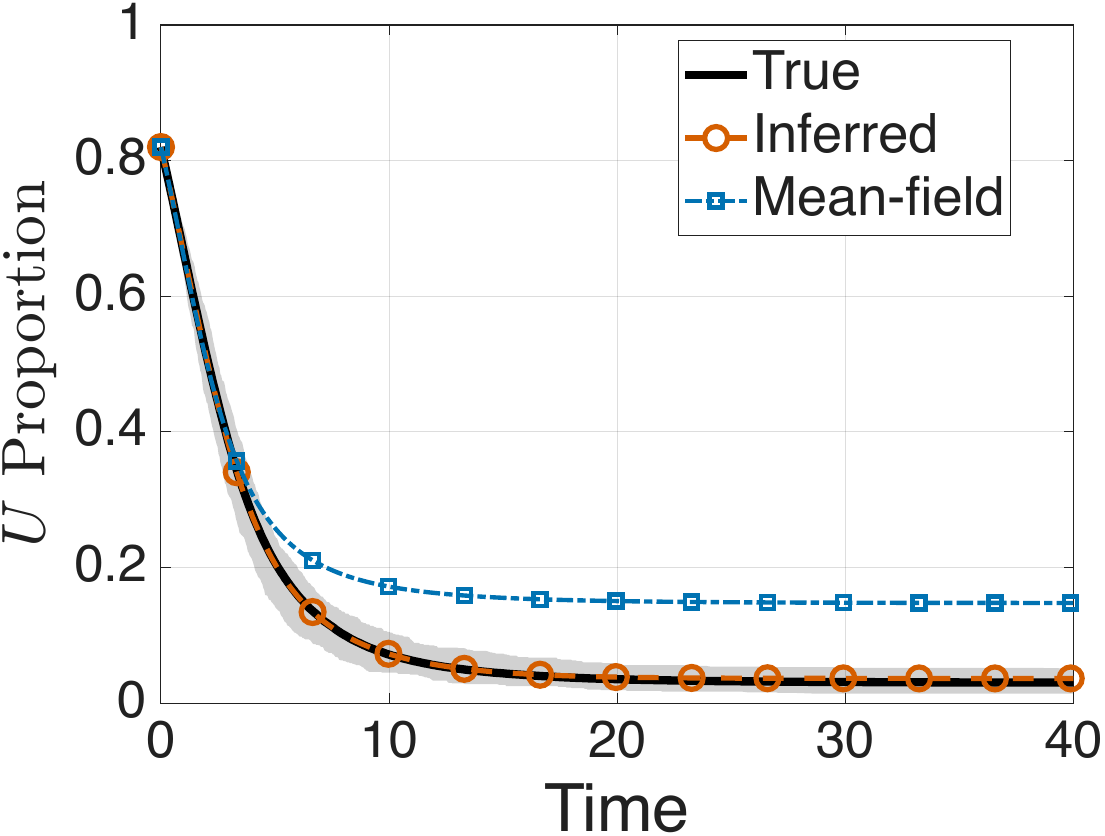}
        \caption{$U$ trajectory}
        \label{subfig:ER_Traj_Compare_U}
    \end{subfigure}
    \hfill
    \begin{subfigure}[t]{0.315\textwidth}
        \centering
        \includegraphics[width=\textwidth]{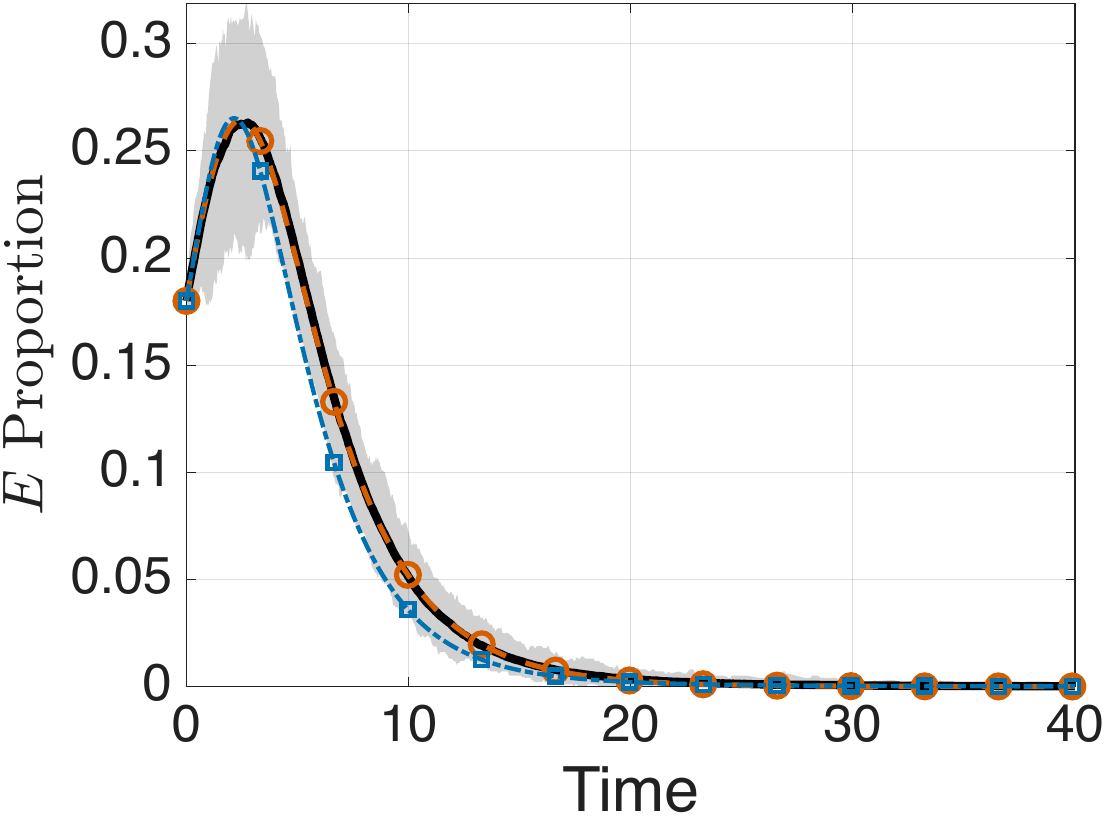}
        \caption{$E$ trajectory}
        \label{subfig:ER_Traj_Compare_E}
    \end{subfigure}
    \hfill
    \begin{subfigure}[t]{0.315\textwidth}
        \centering
        \includegraphics[width=\textwidth]{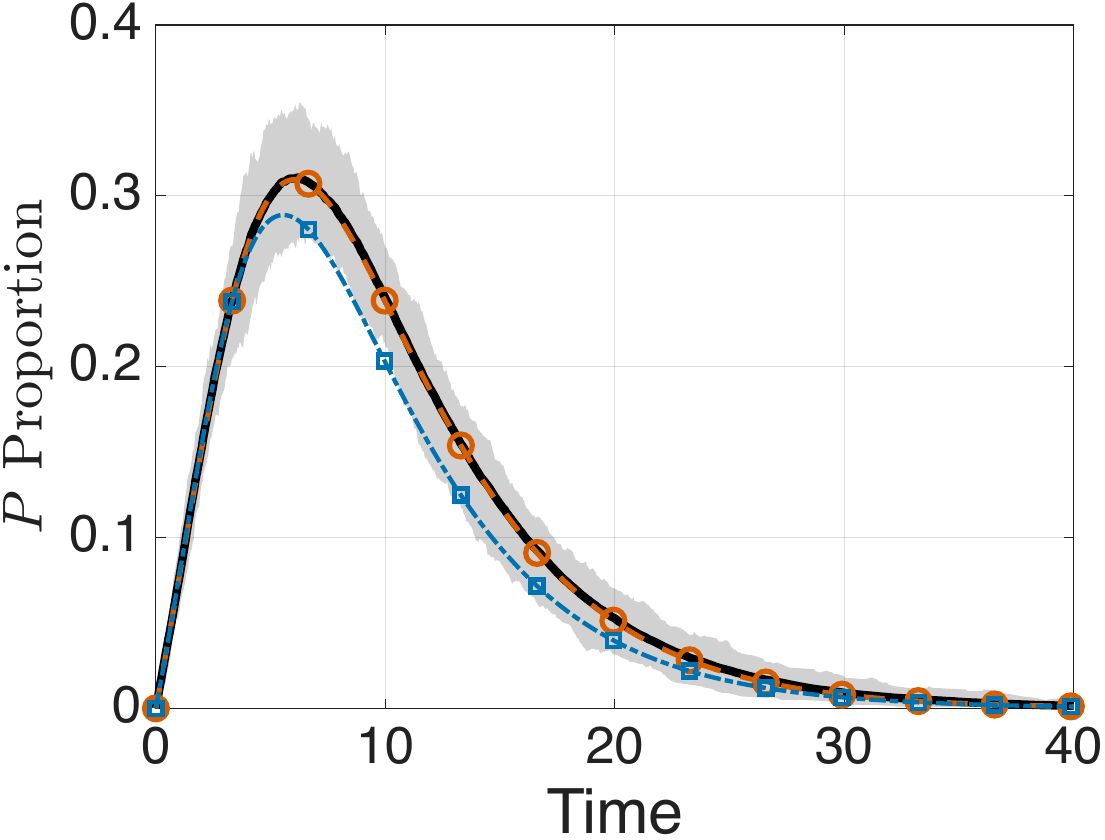}
        \caption{$P$ trajectory}
        \label{subfig:ER_Traj_Compare_P}
    \end{subfigure}
    \caption{Comparison of (a) $U$, (b) $E$, and (c) $P$ trajectories among the sample-averaged stochastic data, the inferred solutions from the learned model based on data, and the solutions from the mean-field approximation model. The average degree is $k = 5$ in this example, with $10$ trajectories used in the training process to obtain the learned model. The same parameters are used for data generation as above. The learned model infers solutions closer to the truth and thus serves as a better approximation model than the mean-field.}
    \label{fig:ER_Traj_Compare}
\end{figure}

In this regard, an advantage of using data-driven modeling techniques is that they yield models that infer solutions that more closely match the observations when network connections are sparse, and the mean-field approximation fails to provide an accurate approximation. Figure \ref{fig:ER_Traj_Compare} shows an example comparing the observed stochastic data (which is considered the truth), the solutions from the mean-field approximation (details in Appendix \ref{app:mean-field}), and the solutions inferred by the learned model. Note that in this example, we used $10$ trajectories on a fixed Erd\H{o}s-R\'{e}nyi network with $k=5$. This scenario corresponds to the extended regime not shown in Figures \ref{fig:ER_FailRate} and \ref{fig:ER_Traj_Err}. The learning in this regime is very stable and accurate, with zero failure and low output errors. As shown in Figure \ref{fig:ER_Traj_Compare}, the learned model is able to predict solutions that are a closer match to the truth across all three compartments, $\bf{U}, \bf{E}$, and $\bf{P}$.

To assess the models learned by WSINDy, Figure \ref{fig:ER_terms_count_weight} shows the counts of identified terms and box plots summarizing the corresponding weights in the equations for $U$, $E$, and $P$, respectively, using data generated from $100$ network realizations with average degree $k=5$. We again use $10$ trajectories, so the learning lands in the stable, considerably accurate regime. Note that each network realization is a different Erd\H{o}s-R\'{e}nyi random network and hence each trajectory is a distinct sample-averaged trajectory over $100$ stochastic realizations. The network heterogeneity, therefore, leads to different learned models. Thus, the variability in Figure~\ref{fig:ER_terms_count_weight} should be interpreted primarily as network-to-network heterogeneity among learned effective models, rather than as estimation uncertainty conditional on a given network. While Figures~\ref{fig:ER_TrajErr} and~\ref{fig:ER_Traj_Compare} provide examples of results on a single fixed network, further analysis distinguishing network-to-network variability from stochastic simulation variability is left for future work. For this reason, we cannot provide a single explicit model, but instead, we here present models as a distribution. 

In Figure \ref{fig:ER_terms_count_weight}, terms appearing in the fully-mixed population model \eqref{model:fully_mixed} for each equation are highlighted by boxes. Note that we use the fully-mixed model as the baseline for comparison because the mean-field model on heterogeneous networks (shown in Appendix~\ref{app:mean-field}) involves additional information, such as the dynamics of the number of nodes in each state with each degree, which is not available to the model learning procedure and also cannot be recovered from population-level data. Thus, model \eqref{model:fully_mixed} provides a more appropriate baseline for interpretability and serves as a sanity check for the learned models.
  
From the weight plots for all three equations, we observe that WSINDy consistently recovers the signs of the terms appearing in the fully-mixed model. This indicates that the corresponding effects, such as the $U \cdot E$ term in the $U$ equation and the $E$ term in the $E$ equation, change the compartmental dynamics in directions consistent with their physical interpretations. Additional terms not present in the fully-mixed model also appear, and those with higher frequencies, such as $E$ in the $U$ equation and $E \cdot P$ in the $E$ equation, suggest that these monomials may serve as effective approximations of network and cross-layer effects. Notably, the learned equation for $P$ aligns closely with the fully-mixed model, with the terms $E$ and $P$ appearing consistently and more dominantly than all other terms, and exhibiting negligible variance across realizations. This aligns with the expectation that network structure exerts a less pronounced influence on offline dynamics, which primarily manifests as an indirect aftereffect of online interactions.
\begingroup
\setlength{\abovecaptionskip}{3pt}
\setlength{\belowcaptionskip}{0pt}

\captionsetup{font=footnotesize, skip=3pt}

\captionsetup[subfigure]{font=scriptsize, skip=0pt}

\begin{figure}[p]
    \centering

    \begin{subfigure}[t]{0.72\textwidth}
        \centering
        \includegraphics[width=\textwidth]{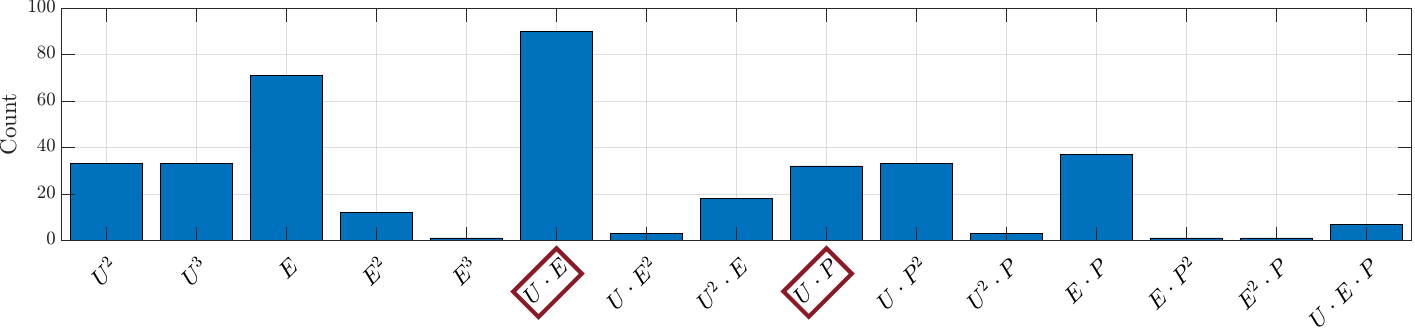}
        \caption{$U$ equation terms count}
        \label{subfig:ER_eq1_terms_count}
    \end{subfigure}\\[1pt]

    \begin{subfigure}[t]{0.72\textwidth}
        \centering
        \includegraphics[width=\textwidth]{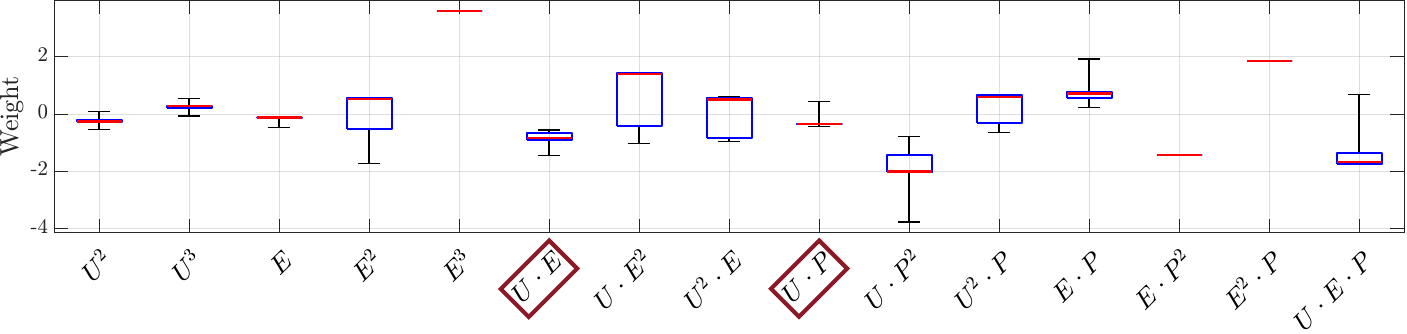}
        \caption{$U$ equation terms weight}
        \label{subfig:ER_eq1_terms_weigth}
    \end{subfigure}\\[1pt]

    \begin{subfigure}[t]{0.72\textwidth}
        \centering
        \includegraphics[width=\textwidth]{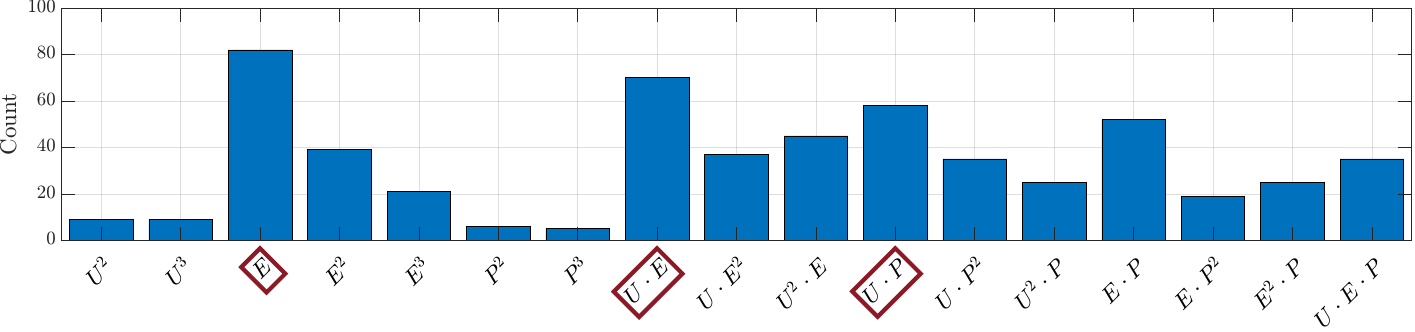}
        \caption{$E$ equation terms count}
        \label{subfig:ER_eq2_terms_count}
    \end{subfigure}\\[1pt]

    \begin{subfigure}[t]{0.72\textwidth}
        \centering
        \includegraphics[width=\textwidth]{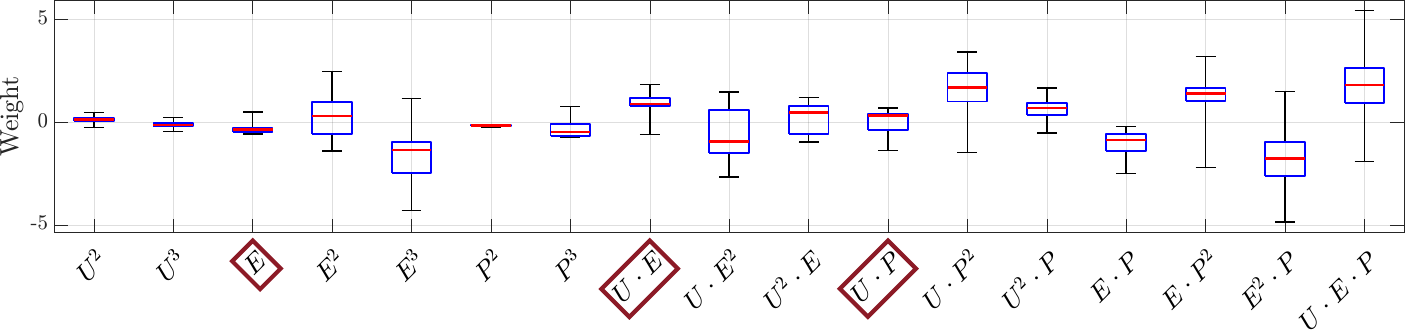}
        \caption{$E$ equation terms weight}
        \label{subfig:ER_eq2_terms_weigth}
    \end{subfigure}\\[1pt]

    \begin{subfigure}[t]{0.72\textwidth}
        \centering
        \includegraphics[width=\textwidth]{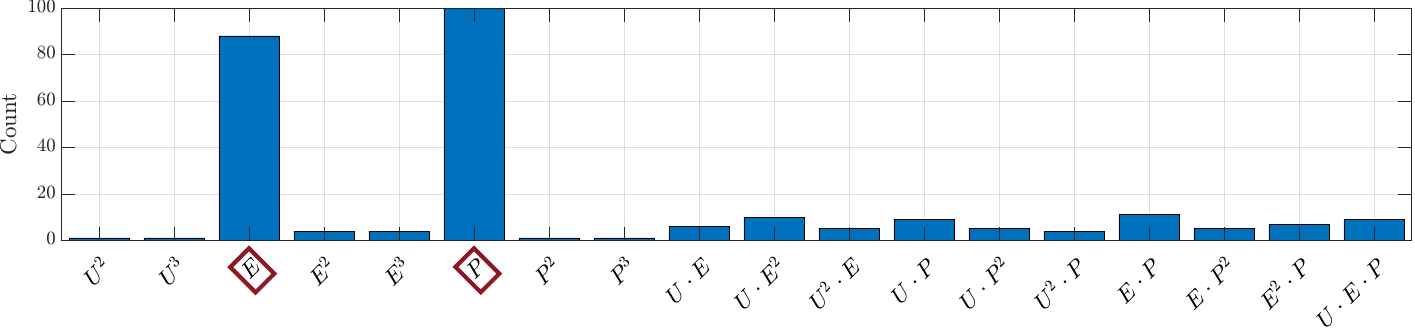}
        \caption{$P$ equation terms count}
        \label{subfig:ER_eq3_terms_count}
    \end{subfigure}\\[1pt]

    \begin{subfigure}[t]{0.72\textwidth}
        \centering
        \includegraphics[width=\textwidth]{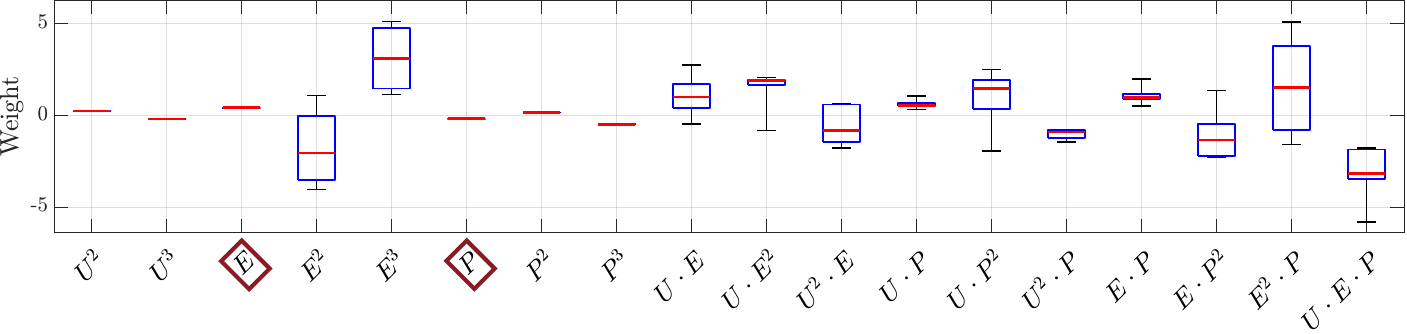}
        \caption{$P$ equation terms weight}
        \label{subfig:ER_eq3_terms_weigth}
    \end{subfigure}

    \caption{Counts of terms appearing across $100$ learned models 
    (panels \ref{subfig:ER_eq1_terms_count}, \ref{subfig:ER_eq2_terms_count}, and \ref{subfig:ER_eq3_terms_count}) 
    and box plots of their corresponding weights 
    (panels \ref{subfig:ER_eq1_terms_weigth}, \ref{subfig:ER_eq2_terms_weigth}, and \ref{subfig:ER_eq3_terms_weigth}). Each learned model is identified from data generated on one of $100$ independent network realizations with average degree $k = 5$, where each trajectory is averaged over $100$ stochastic realizations on that network. The spread in the box plots therefore primarily reflects variability across network realizations. Other parameters used for data generation are the same as those described above. The boxed terms correspond to those appearing in the fully-mixed population model \eqref{model:fully_mixed} for the respective equations. The signs of these learned terms align with those in the model. In particular, the learned equation for $P$ aligns especially well, with its corresponding terms appearing dominantly relative to others. Additional terms that appear with high frequency in learning suggest important terms to include in approximations that capture network effects.}
    \label{fig:ER_terms_count_weight}
\end{figure}
\endgroup

With the goal to learn a model that better matches the data, in practice, we may exploit the learning results, for example, as in Figure \ref{fig:ER_terms_count_weight}, to build a more restricted library for WSINDy by selecting terms that appear with frequencies above a certain threshold, indicating such terms are likely to carry more importance in describing the dynamics. Specific to the example in Figure \ref{fig:ER_terms_count_weight}, upon investigation, we found that using a threshold of $30$, selecting library terms with counts above this value, yields the best learning performance in terms of minimizing the output errors. In these more restricted runs, some learned terms begin to appear more deterministically, appearing $100$ times across $100$ learned models. We observe that these terms are $U$, $E$, $U^2$ and $U^3$ in the $U$ equation, $E$ and $U \cdot E$ in the $E$ equation, and $E$ and $P$ in the $P$ equation. This suggests that these terms should be expected in the approximation model for the stochastic dynamics on networks that we are considering.

\section{Conclusion} \label{sec:conclusion}

In this work, we present a systematic study of using the data-driven system identification method, WSINDy, to learn the underlying equations governing coupled online-offline social dynamics. We specifically examine two models: a continuum model consisting of a system of ODEs for a fully-mixed population, and a stochastic model in which individuals exhibit hybrid online-offline states and interact through networks in the online domain. 

The fully-mixed model provides a testbed in which the ground-truth model is known. Thus, we examine the TPR and parameter error, in addition to equation and output errors, as performance metrics, and present studies on the influence of noise. We observe that adding a small number of additional trajectories to the training set can significantly improve learning performance with greater noise tolerance across all metrics. Additionally, contrary to intuition, adding more trajectories beyond the minimal number of 3 does not yield further improvement. 

For the stochastic model, we specifically focus on Erd\H{o}s-R\'{e}nyi networks and use output error as the measurement for performance. We find the same takeaways from these experiments: a minimal number of trajectories can improve learning, meaning producing predicted solutions that are both physically meaningful and numerically solvable and that present a closer match to the true trajectories, as measured by output error. However, again, little additional benefit is observed beyond 3 trajectories. Another advantage of data-driven modeling is that, in the sparse network scenario, which is often true in empirical settings, it yields better approximations that better predict empirical dynamics than analytically derived mean-field models.

Our work provides a first effort to connect social dynamics data on network structure to efficient models, aiming to provide an empirical means to reveal the effective governing mechanisms underlying real-world dynamics. The fully-mixed model, although simple, already exhibits variability in the system identification performance in the presence of noise. One limitation is that although our learning procedure mitigates the most direct mass-conservation issue by using only the active trajectories in the regression, noisy observations may still fall slightly outside their physically meaningful ranges. While our work presented computational analysis, one potential future endeavor is to use the zero-noise case of this system, or other mean-field approximation models, as a diagnostic to study the theoretical guarantees of learning performance in relation to properties such as the rank, data sufficiency, and number of trajectories in the data. Our results also suggest several other useful directions for future research, including developing an efficient method for selecting which trajectory information to use in model learning and applying our method to real-world network dynamics data with missing information or hidden states. 

\section*{Code availability}
The code used to generate the results in this paper is available at
\url{https://github.com/Moyi-Tian/WSINDy-NetworkDynamics}.

\appendix



\section{Assessing identifiability} \label{app:eq_criterion}

For the fully-mixed model, we investigate the \textit{structural identifiability} of the system with $(U,E,D,P,R), (U,E,P,R)$, and $(E, P)$, respectively. The Julia \textit{StructuralIdentifiability} package indicates that all of these full/reduced models are globally identifiable with respect to all terms and parameters. This structural identifiability result is an idealized property of the model and observation structure, and does not by itself guarantee accurate recovery from finite, noisy data using a particular estimator.

However, in realistic settings, recovery also depends on several factors, including noise, the amount of measured data, and the estimation method. Inspired by the $(\mathbf{e},\mathbf{q})$-identifiability criterion \cite{Heitzmanbreen2025}, which assesses \textit{practical identifiability} of a system given the available data and parameter estimation method, we apply a slightly modified version of this criterion to our setting. In our case, models are obtained through system identification rather than direct parameter estimation. Such a modification is necessary because the learned model may assign nonzero weights to terms whose true coefficients are zero, which would lead to division-by-zero issues if the criterion were applied to each parameter individually.

We therefore evaluate identifiability using the relative $\ell^2$ error between the vectorized learned and true coefficient matrices:
\[
\frac{1}{r}\sum_{i=1}^{r}
\frac{\left\|\operatorname{vec}\left(\hat{\mathbf{W}}_i(\mathbf{e})\right)-\operatorname{vec}\left(\mathbf{W}\right)\right\|_2}{\left\|\operatorname{vec}\left(\mathbf{W} \right)\right\|_2} < \mathbf{q},
\]
where $\mathbf{e}$ denotes the additive error ratio, $\operatorname{vec}\left(\hat{\mathbf{W}}_i(\mathbf{e})\right)$ is the vectorization of the estimated coefficient matrix from the $i$-th realization, $\operatorname{vec}\left(\mathbf{W} \right)$ is the vectorization of the true coefficient matrix, and $r$ is the number of realizations. 

This metric is interpreted as follows: at noise level $\mathbf{e}$, the average relative $\ell^2$ error of the vectorized estimated coefficient matrix is bounded above by $\mathbf{q}$. Thus, unlike structural identifiability, this modified criterion should be interpreted as a finite-data estimation-accuracy diagnostic for the chosen WSINDy procedure, rather than as a sufficient condition for identifiability. In real data applications, the true coefficient matrix is unknown, so an analogous diagnostic would require replacing $\mathbf{W}$ by a reference estimate, such as the empirical mean of the estimated coefficient matrices across repeated fits.

\begin{figure}[!htbp]
    \centering
    \begin{subfigure}[t]{0.33\textwidth}
        \centering
        \includegraphics[width=\textwidth]{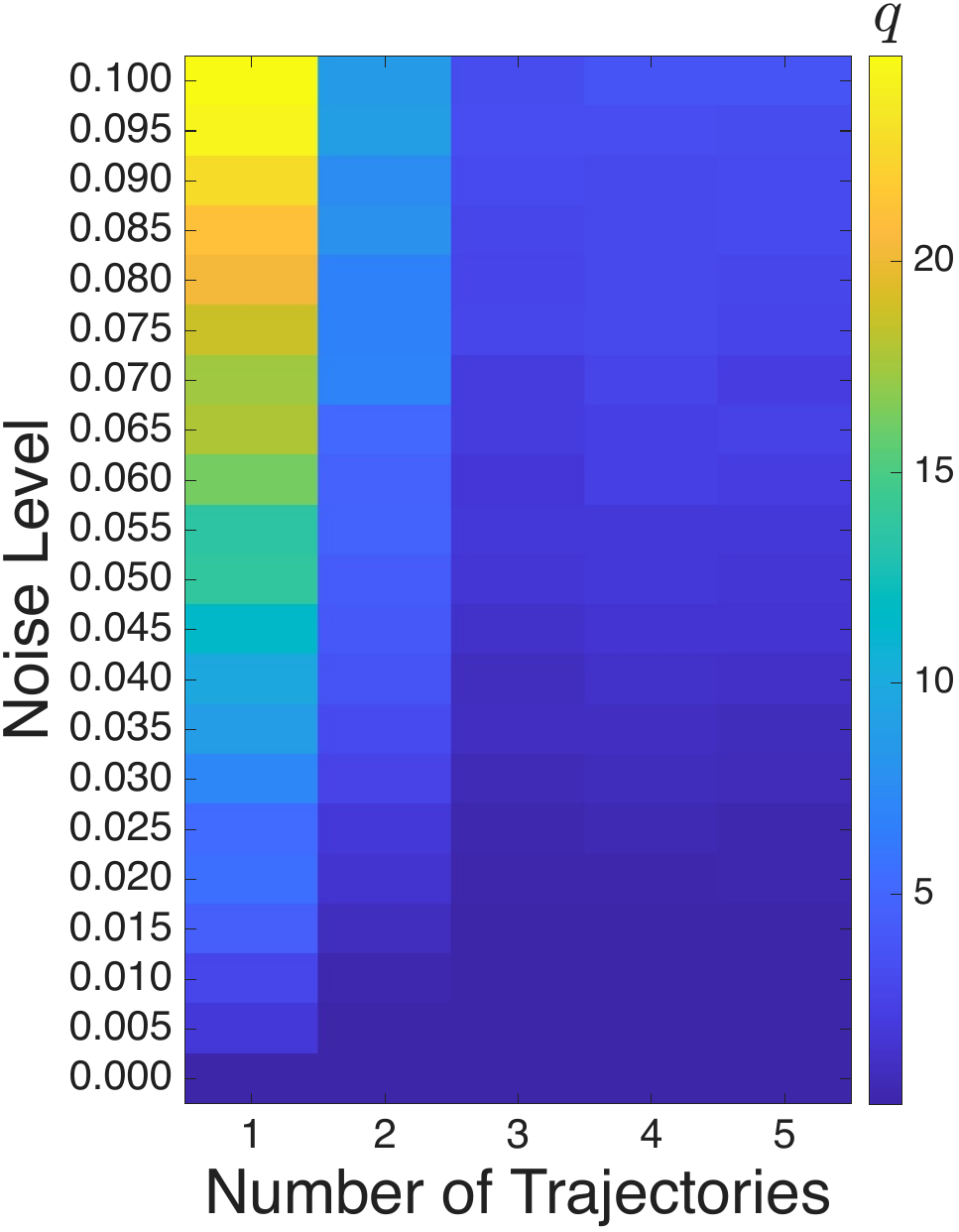}
        \caption{True coefficients}
        \label{subfig:E_q_surface_Noise_NumTrajs}
    \end{subfigure}
    \hspace{50pt}
    \begin{subfigure}[t]{0.31\textwidth}
        \centering
        \includegraphics[width=\textwidth]{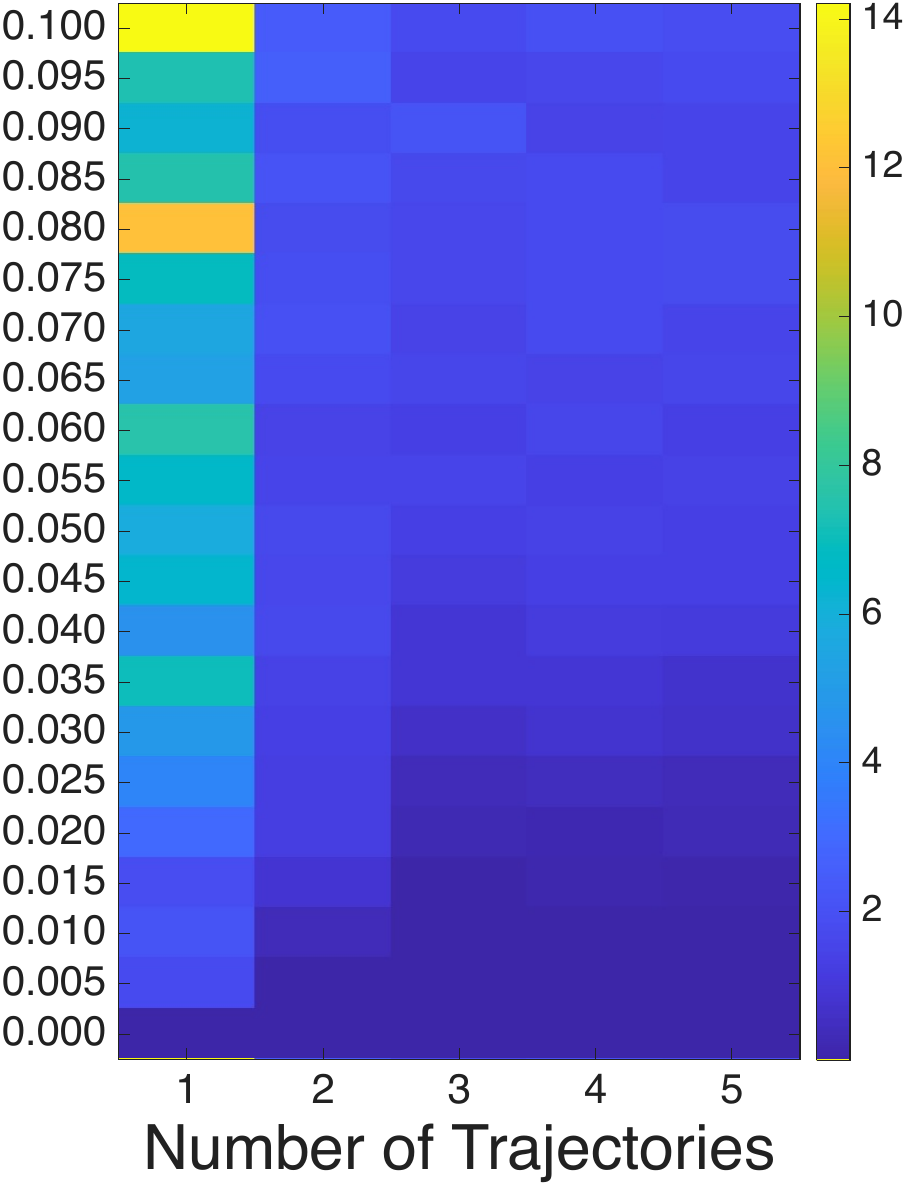}
        \caption{Empirical mean}
        \label{subfig:E_q_surface_empMean}
    \end{subfigure}
    \caption{$(\mathbf{e},\mathbf{q})$-identifiability summary for $100$ realizations with varying initial conditions in $E(0)$. (a) Heatmap of the average relative $\ell^2$ error of the vectorized learned coefficient matrix, using the true coefficient matrix as the reference. (b) Analogous heatmap using the empirical mean of the learned coefficient matrices as the reference. Both panels are shown over noise level and number of training trajectories. The same parameters are used for data generation as in Figure~\ref{fig:E_TPS_ParamErr_Noise_NumTrajs}. Both panels show that $q$ decreases substantially as the number of trajectories increases up to $3$, with the most pronounced drop from $1$ to $2$ and comparatively small changes beyond $3$.}
    \label{fig:E_q_surface_CrossSect}
\end{figure}

Figure~\ref{subfig:E_q_surface_Noise_NumTrajs} shows the $q$-surface heatmap associated with the $(\mathbf{e},\mathbf{q})$-identifiability criterion using the true coefficient matrix as the reference. Figure~\ref{subfig:E_q_surface_empMean} shows the analogous empirical-reference heatmap, where the true coefficient matrix is replaced by the empirical mean of the learned coefficient matrices at each noise level and number of trajectories. The results indicate a substantial decrease in $q$ when the number of trajectories increases from $1$ to $3$.

\section{Mean-field model on heterogeneous networks} \label{app:mean-field}

We present the mean-field approximation for the stochastic model on heterogeneous networks, which is used to compare with the data-driven model in Section \ref{subsec:results_stoch}. We refer to the details of the derivation in \cite{Tian2025}.

For each compartment label $\mathbf{X}\in\{\mathbf{U},\mathbf{E},\mathbf{D},\mathbf{P},\mathbf{R}\}$, let $X$ be the expected population fraction in compartment $\mathbf{X}$, let $X_n$ be the fraction of people in compartment $\mathbf{X}$ with degree $n$, and let $\hat{n}$ be the fraction of people with degree $n$. In the numerical comparison with Erd\H{o}s-R\'{e}nyi networks, the degree fractions $\hat n$ are computed directly from the realized network. Suppose that the maximum degree of any node in the network is $K_{\text{max}}$. For $n=1,\ldots, K_{\text{max}}$, the mean-field approximation is
\begin{align*}
\left\{\begin{array}{ll}
    \frac{dU_n}{dt}&=-\tau n U_n \pi_E - \theta U_n P,\\[4pt]
    \frac{dE_n}{dt}&=\tau n U_n \pi_E + \theta U_n P - (\eta+\gamma_i) E_n,\\[4pt]
    \frac{dD_n}{dt}&=(\eta+\gamma_i) E_n,\\[4pt]
    \frac{dP(t)}{dt}& = \eta E-\gamma_p P,\\[4pt]
    \frac{dR(t)}{dt}& = \gamma_p P,
\end{array}\right.
\end{align*}
where $\pi_E$ is the probability that a randomly chosen stub is connected to a stub of an online engaged individual, and is equal to
\[
\pi_E = \frac{\sum_{j=1}^{K_{\text{max}}} j E_j}{\sum_{j=1}^{K_{\text{max}}} j \hat{j}}.
\]

\section{Variability in fully mixed WSINDy performance}
\label{app:variability_fullymixed}

\begin{figure}[H]
    \centering
    \begin{subfigure}[t]{0.38\textwidth}
        \centering
        \includegraphics[width=\textwidth]{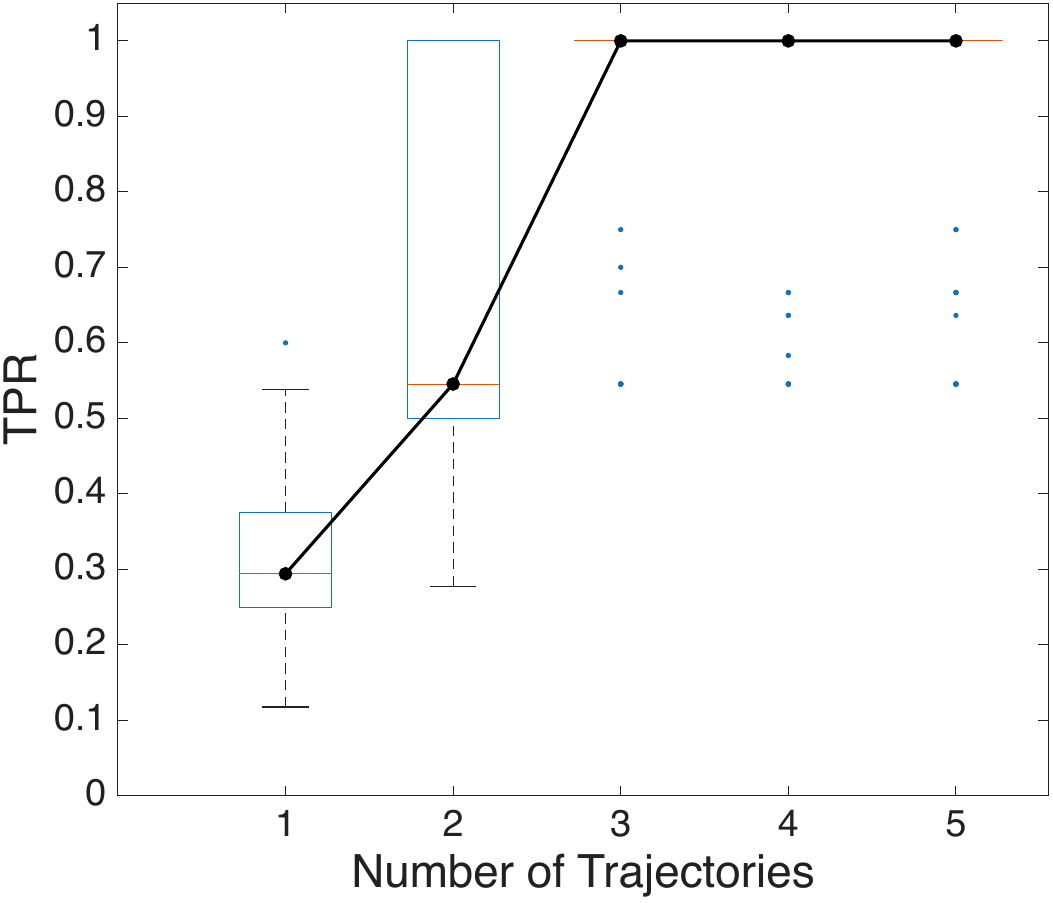}
        \caption{TPR}
        \label{subfig:app_TPR_CrossSect}
    \end{subfigure}
    \hspace{40pt}
    \begin{subfigure}[t]{0.39\textwidth}
        \centering
        \includegraphics[width=\textwidth]{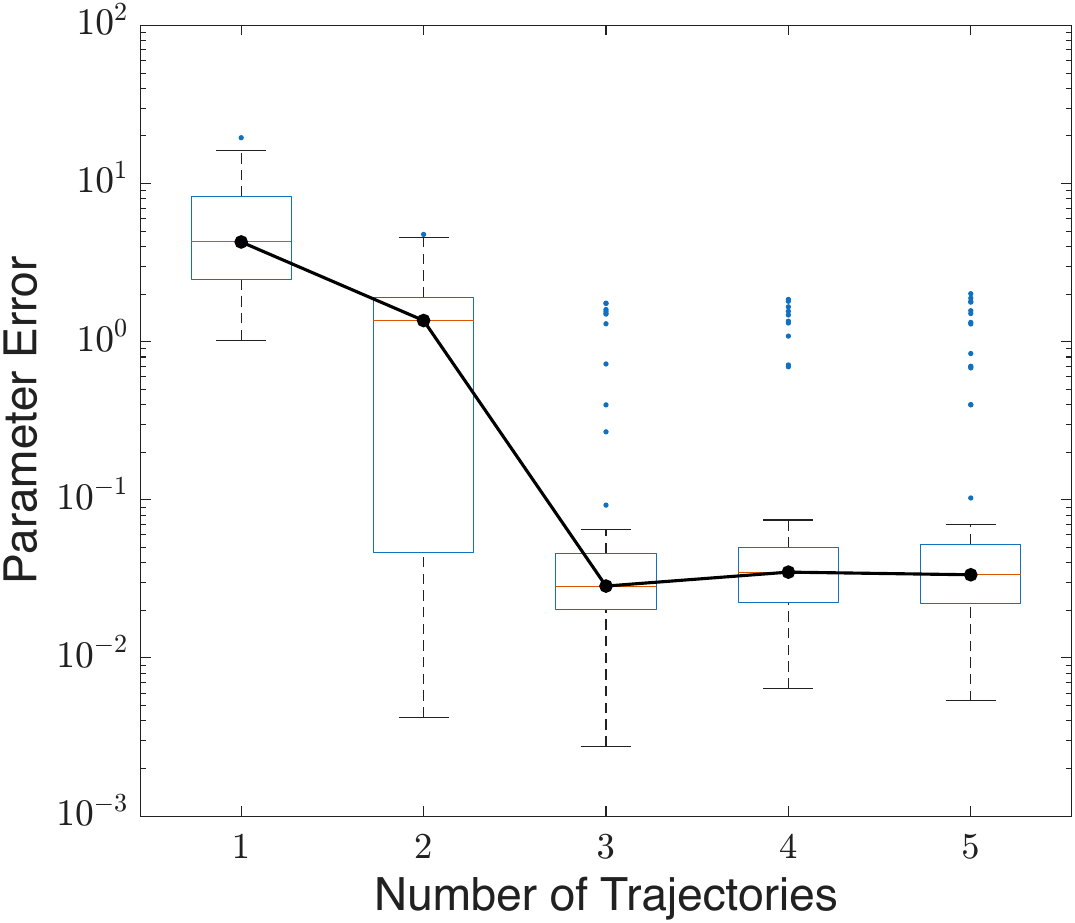}
        \caption{Parameter error}
        \label{subfig:app_ParamErr_CrossSect}
    \end{subfigure}
    \caption{Cross-sectional distributions of (a) TPR and (b) parameter error over $100$ independent noise realizations at noise level $0.02$, corresponding to the heatmaps in Figure~\ref{fig:E_TPS_ParamErr_Noise_NumTrajs}. The distributions show the main improvement from adding the first few trajectories and diminishing changes after $3$ trajectories. For TPR, some boxes collapse to a line at $\mathrm{TPR}=1$ because the interquartile range is zero.}
    \label{fig:TPR_ParamErr_CrossSect}
\end{figure}

\begin{figure}[H]
    \centering
    \begin{subfigure}[t]{0.38\textwidth}
        \centering
        \includegraphics[width=\textwidth]{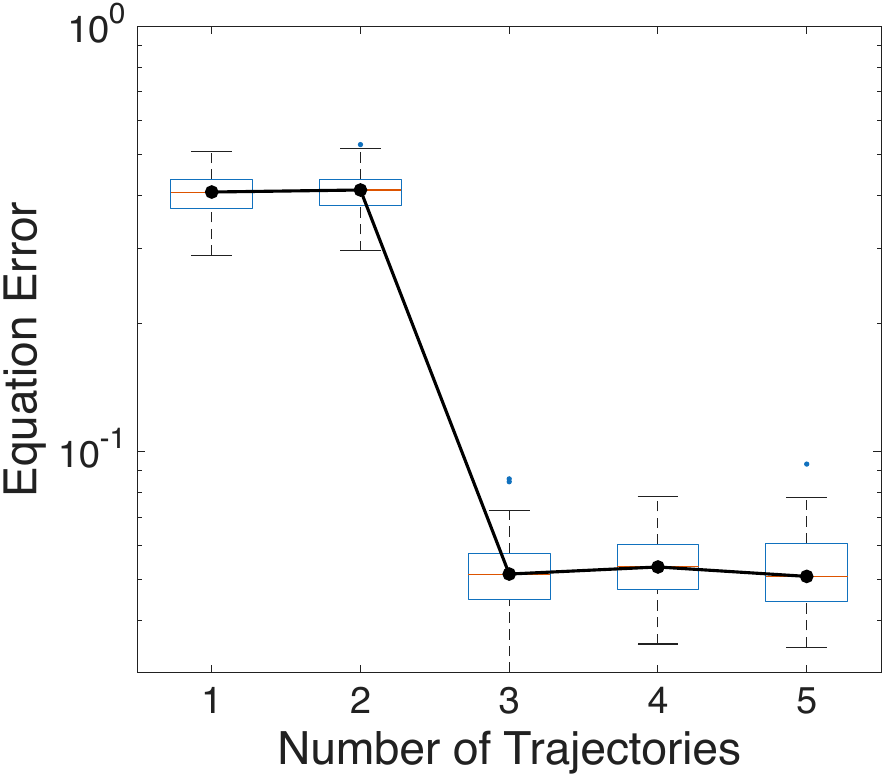}
        \caption{$E$ equation}
        \label{subfig:EquaErr_CrossSect_E}
    \end{subfigure}
    \hspace{50pt}
    \begin{subfigure}[t]{0.35\textwidth}
        \centering
        \includegraphics[width=\textwidth]{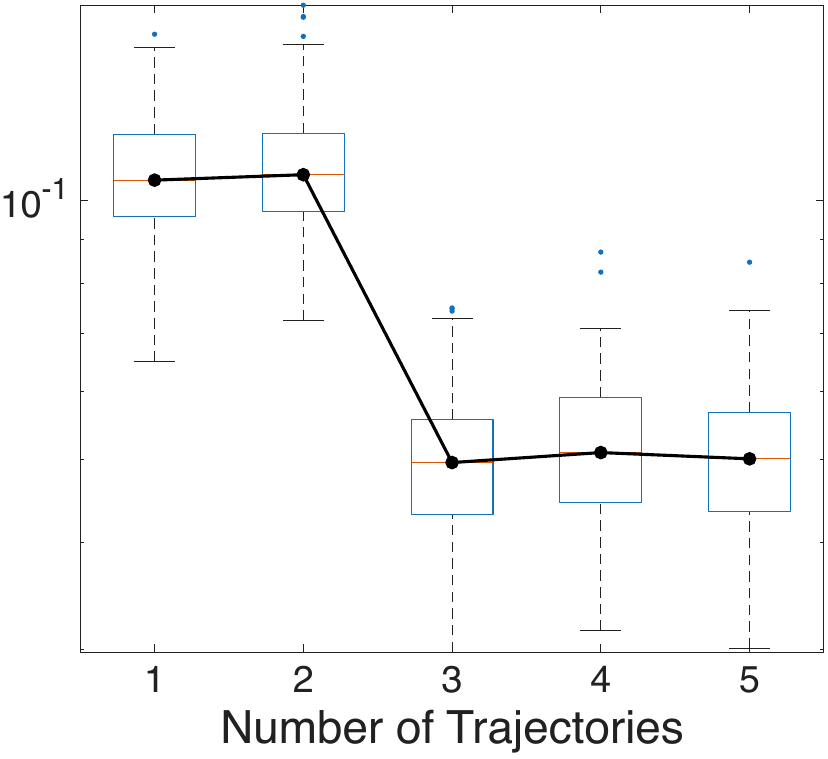}
        \caption{$P$ equation}
        \label{subfig:EquaErr_CrossSect_P}
    \end{subfigure}
    \caption{Cross-sectional distributions of equation error for (a) $E$ and (b) $P$ over $100$ independent noise realizations at noise level $0.05$, corresponding to the heatmaps in Figure~\ref{fig:E_EquaErr_Noise_NumTrajs}. These cross-sections show the main reduction in error occurs when increasing from $2$ to $3$ trajectories, with comparatively small changes beyond that.}
    \label{fig:EquaErr_CrossSect}
\end{figure}

\begin{figure}[H]
    \centering
    \begin{subfigure}[t]{0.38\textwidth}
        \centering
        \includegraphics[width=\textwidth]{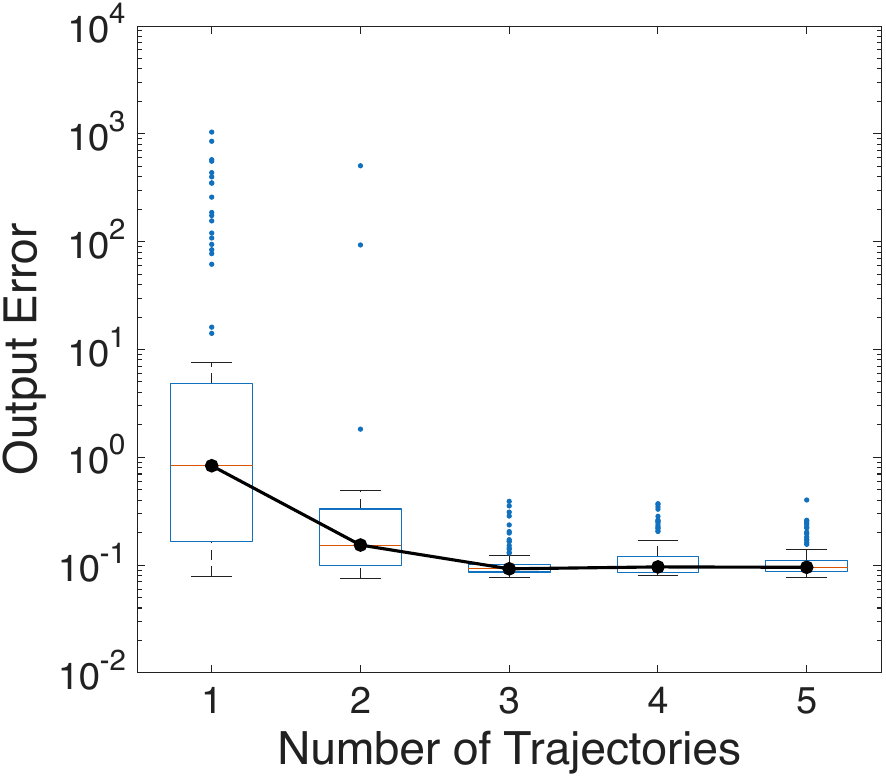}
        \caption{$E$ trajectory}
        \label{subfig:OutErr_CrossSect_E}
    \end{subfigure}
    \hspace{50pt}
    \begin{subfigure}[t]{0.35\textwidth}
        \centering
        \includegraphics[width=\textwidth]{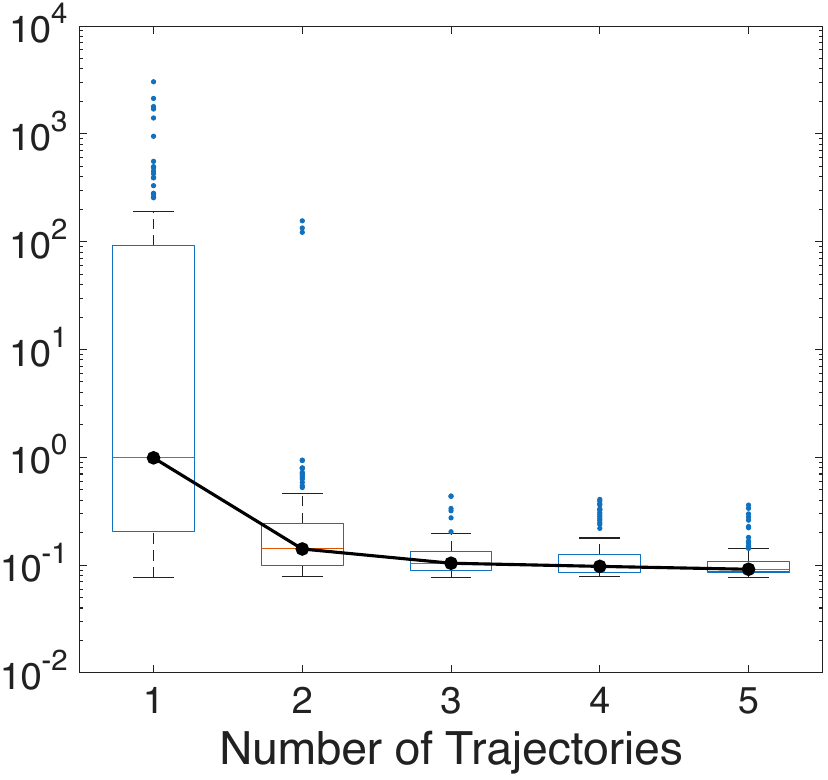}
        \caption{$P$ trajectory}
        \label{subfig:OutErr_CrossSect_P}
    \end{subfigure}
    \caption{Cross-sectional distributions of output error for (a) $E$ and (b) $P$ over $100$ independent noise realizations at noise level $0.08$, corresponding to the heatmaps in Figure~\ref{fig:E_OutErr_Noise_NumTrajs}. These cross-sections support the observation that improvements are most pronounced from $1$ to $2$ and then from $2$ to $3$ trajectories, with limited additional changes afterward.}
    \label{fig:OutErr_CrossSect}
\end{figure}

The main results in Figures~\ref{fig:E_TPS_ParamErr_Noise_NumTrajs}--\ref{fig:E_OutErr_Noise_NumTrajs} report median performance over $100$ noise realizations. To illustrate the variability behind these summaries, Figures~\ref{fig:TPR_ParamErr_CrossSect}--\ref{fig:OutErr_CrossSect} show representative cross-sectional boxplots at selected noise levels. These plots are intended as descriptive variability summaries rather than exhaustive statistical tests over all noise-level and trajectory-count pairs. Boxes show the interquartile range, center lines show the median, whiskers extend to the most extreme points within $1.5$ times the interquartile range, and points outside this range are shown individually. Black markers connect the median values across the number of training trajectories. In particular, these selected cross-sections provide additional evidence that the most pronounced improvements occur when the number of trajectories is increased from $1$ to $2$ or from $2$ to $3$, while changes beyond $3$ trajectories are comparatively small.

\section{Variability in stochastic network WSINDy performance}
\label{app:variability_stochastic}

\begin{figure}[H]
    \centering
    \begin{subfigure}[t]{0.325\textwidth}
        \centering
        \includegraphics[width=\textwidth]{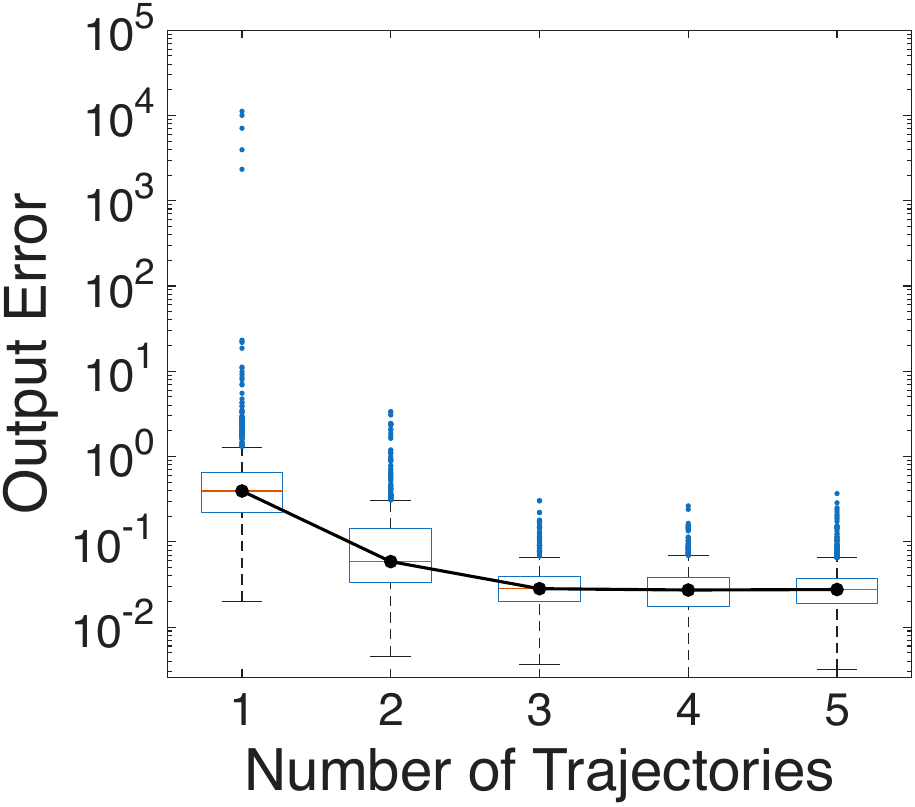}
        \caption{$U$ trajectory}
        \label{subfig:app_ER_TrajErr_U}
    \end{subfigure}
    \hfill
    \begin{subfigure}[t]{0.3\textwidth}
        \centering
        \includegraphics[width=\textwidth]{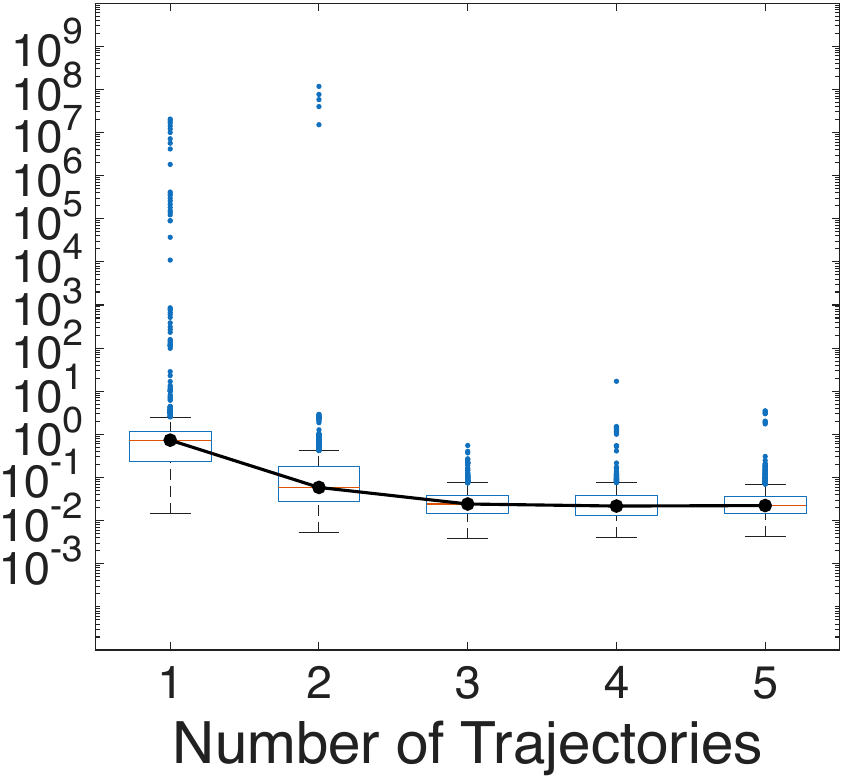}
        \caption{$E$ trajectory}
        \label{subfig:app_ER_TrajErr_E}
    \end{subfigure}
    \hfill
    \begin{subfigure}[t]{0.3\textwidth}
        \centering
        \includegraphics[width=\textwidth]{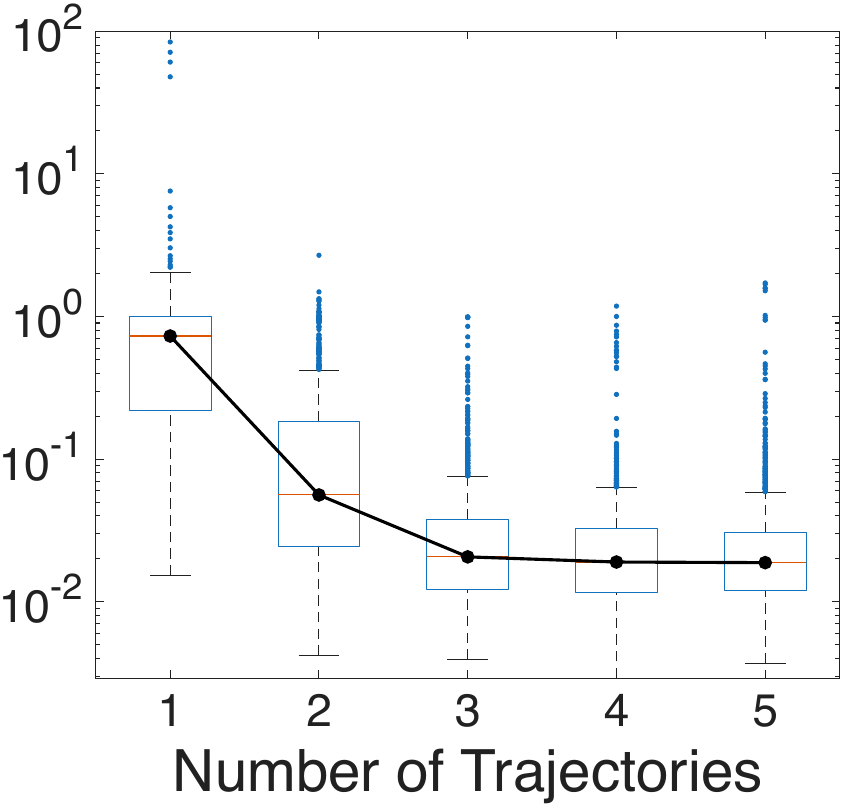}
        \caption{$P$ trajectory}
        \label{subfig:app_ER_TrajErr_P}
    \end{subfigure}
    \caption{Cross-sectional distributions of output error for (a) $U$, (b) $E$, and (c) $P$ at average degree $k=5$, corresponding to the heatmaps in Figure~\ref{fig:ER_Traj_Err}. The plots show variability across random network realizations and test sample-averaged trajectories while preserving the same qualitative trend as the mean heatmaps. The main improvement occurs by $3$ trajectories, with smaller and less systematic changes beyond that.}
    \label{fig:app_ER_TrajErr_CrossSect}
\end{figure}

The heatmaps in Figure~\ref{fig:ER_Traj_Err} report mean output errors over random network realizations and test datasets. To show the variability behind these mean summaries, Figure~\ref{fig:app_ER_TrajErr_CrossSect} presents representative cross-sectional boxplots at average degree $k=5$. This value is chosen as a representative sparse-network case and is also used in the trajectory comparison (Figure~\ref{fig:ER_Traj_Compare}) and learned-term analyses (Figure~\ref{fig:ER_terms_count_weight}) in Section~\ref{subsec:results_stoch}. This selected cross-section is used as representative evidence for the trend observed in the heatmaps: output errors decrease most noticeably up to $3$ training trajectories, while changes beyond $3$ trajectories are smaller and less systematic.

\section*{Acknowledgments}
This research was supported by the NIH-NIGMS Division of Biophysics, Biomedical Technology and Computational Biosciences [grant R35GM149335]; the U.S. National Science Foundation Division of Mathematical Sciences [grant 2042413]; and the Air Force Office of Scientific Research Multidisciplinary University Research Initiative [grant FA9550-22-1-0380].

We thank our colleagues, Nora Heitzman-Breen, Rainey Lyons, and April Tran, for valuable discussions and feedback at various stages of this work.

\FloatBarrier
\bibliographystyle{plain}
\bibliography{WDynOnNet}

\end{document}